\documentclass[12pt] {article} 

\newtheorem{theorem}{Theorem}[section]
\newtheorem{lemma}{Lemma}[section]

\newcommand{\eqnsection}{
   \renewcommand{\theequation}{\thesection.\arabic{equation}}
   \makeatletter
   \csname @addtoreset\endcsname{equation}{section}
   \makeatother}
 
\def \ov{\overline}

\def \be{\begin{equation}}
\def \ee{\end{equation}}
\def \bt{\begin{theorem}} 
\def \et{\end{theorem}}
\def \bl{\begin{lemma}} 
\def \el{\end{lemma}}
\def \bea{\begin{eqnarray}}
\def \eea{\end{eqnarray}}
\def \bas{\begin{eqnarray*}}
\def \eas{\end{eqnarray*}}

\def \al{\alpha}
\def \bb{\beta}

\def \de{\delta}
\def \De{\Delta}
\def \ep{\epsilon}

\def \la{\lambda}

\def \si{\sigma}

\def \ze{\zeta}

\def \ff{\infty}
\def \wh{\widehat}
\def \wt{\widetilde}

\def \rar{\rightarrow}

\def \cd{\,\cdot\,}

\def \BB{{\cal B}}

\def \II{{\cal I}}

\def \LL{{\cal L}}

\def \PP{{\cal P}}

\def \TT{{\cal T}}

\def \E{E}

\def \({\left(}
\def \){\right)}

\def \lc{\left\{}
\def \rc{\right\}}

\def \nn{\nonumber}
\def \Proof{\noindent{\bf Proof $\,$ }}

\def \bc{\begin{center} }
\def \ec{\end{center} }
\def \bs{\begin{slide} }
\def \es{\end{slide} }

\def\square{{\vcenter{\vbox{\hrule height.3pt
        \hbox{\vrule width.3pt height5pt \kern5pt
           \vrule width.3pt}
        \hrule height.3pt}}}}
\def\qed{{\hfill $\square$ \bigskip}}

\eqnsection
\begin{document}

\def\wh{\widehat}
\def\ol{\overline}

\title{A CLT  for the $L^{2}$ modulus of continuity \\of  Brownian local time}

\author{   Xia Chen,    Wenbo Li,  Michael B. Marcus and  Jay Rosen\thanks
     {  The research of all four authors was  supported, in part, by grants from the National Science
Foundation. The research of Professors Marcus and Rosen was also  supported, in part, by grants from PSC-CUNY.}}


\maketitle

\bibliographystyle{amsplain}

\begin{abstract} { Let $\{L^{ x }_{ t}\,;\,(x,t)\in R^{ 1}\times  R^{  1}_{ +}\}$ denote the local time of  Brownian motion and  	
\[
\alpha_{ t}:=\int_{-\infty}^{\infty} ( L^{ x}_{ t})^{ 2}\,dx .
\]
Let $\eta=N(0,1)$ be  
independent of $\alpha_{ t}$. For each fixed $t$
\[  { \int_{-\infty}^{\infty} ( L^{ x+h}_{t}- L^{ x}_{ t})^{ 2}\,dx- 4ht\over h^{ 3/2}}
\stackrel{\mathcal{L}}{\rightarrow}\( {64 \over  3} \)^{ 1/2}\sqrt{\alpha_{t}}\,\,\eta,
\]
 as $h\rar 0$.  
 Equivalently 
\[ { \int_{-\infty}^{\infty} ( L^{ x+1}_{t}- L^{ x}_{ t})^{ 2}\,dx- 4t\over t^{ 3/4}}
\stackrel{\mathcal{L}}{\rightarrow}\( {64 \over  3} \)^{ 1/2}\sqrt{\alpha_{1}}\,\,\eta ,
\]
as $t\rar\infty$.}

\end{abstract}

\section{Introduction}

 In   \cite{lp}  almost sure limits are obtained for the 
$L^{ p}$ moduli  of continuity of  local times of a very wide class of
symmetric L\'evy processes. For Brownian motion the result is as follows:  Let $\{L^{ x }_{ t}\,;\,(x,t)\in R^{ 1}\times  R^{  1}_{ +}\}$ denote the local time of  Brownian motion. Then    for all 
$  p\ge 1$, and all
$t\in R_+$
\be
\lim_{ h\downarrow 0}  \int_{a}^{ b} \bigg|{  L^{ x+h}_{ t} -L^{ x }_{
t}\over\sqrt{h}}\bigg|^p\,dx =2^pE(|\eta|^p)
\int_a^b |L^{ x }_{ t}|^{ p/2}\,dx\label{as.1}
\ee for all
$a,b
$ in the extended real line  almost surely, and also  in $L^m$, $m\ge 1$.
(Here $\eta$  is a   normal random variable with mean zero and variance
one.) When $p=2$ and $a=-\ff,\,b=\ff$ we can write (\ref{as.1}) in the form
\begin{equation}
\lim_{ h\downarrow 0}  \int_{-\ff}^{\ff}  {  (L^{ x+h}_{ t} -L^{ x }_{
t})^{ 2}\over h}\,dx =4t  \hspace{.2 in}\mbox{ a. s.} \label{rp3.1}
\end{equation}
  This result in (\ref{as.1}) uses the
Eisenbaum Isomorphism Theorem, see e.g. \cite[Theorem 8.1.1]{book}, and is a consequence of a  similar result for the
Ornstein--Uhlenbeck process, (the stationary Gaussian process $\{G(x),x\in
R^1\}$,   with
$E(G(x)-G(y))^2=2( 1-e^{ -|x-y|})$), which is that    for all $p\ge 1$
\be\quad
\lim_{h\to 0} \int_a^b\bigg|\frac{G(x+h)-G(x)}{\sqrt h}\bigg|^p\,dx=E|\eta |^p (b-a)\qquad \forall \,a,b\in R^{1}\, \mbox{  a.s.}\label{1.3}
\ee  
This is also obtained in \cite{lp},   in which this question is  considered   for a very large class of Gaussian processes. The right--hand side of  (\ref{1.3}) is the expected value of the left--hand side.  Consequently, (\ref{1.3}) can be thought of as a law of large numbers. 
  In   \cite{clt}  we   consider the central limit theorem for the  left--hand side of  (\ref{1.3}).   For the Ornstein--Uhlenbeck, when  $p=2$, we get 
\begin{equation}
\lim_{ h\downarrow 0} { \int_{ a}^{ b} (G( x+h)- G( x))^{ 2}\,dx-2h(
b-a)\over h^{ 3/2}}\stackrel{\mathcal{L}}{=}(16/3)^{
1/2}(b-a)\,\eta\label{7.1}.
\end{equation}

   The argument involving the Eisenbaum Isomorphism Theorem that is used in \cite{lp} to show that (\ref{1.3}) implies (\ref{as.1}) does not work to show that (\ref{7.1}) implies a similar result for the local times of Brownian motion. In this paper we obtain a central limit theorem corresponding to (\ref{as.1}) by considering moments of
   \begin{equation}
 \int ( L^{ x+1}_{t}- L^{ x}_{ t})^{ 2}\,dx.
   \end{equation} 
   (An integral sign without limits is to be read as $\int_{-\ff}^{\ff}$.)
   
 Let
\begin{equation}
\al_{ t}=\int ( L^{ x}_{ t})^{ 2}\,dx\label{5.1w}
\end{equation} and let $\eta=N(0,1)$ be  
independent of $\al_{ t}$. We have the following weak convergence  results.

\begin{theorem}\label{theo-clt2} For each fixed $t$
\begin{equation} { \int ( L^{ x+h}_{t}- L^{ x}_{ t})^{ 2}\,dx- 4ht\over h^{ 3/2}}
\stackrel{\mathcal{L}}{\rightarrow}c\sqrt{\al_{t}}\,\,\eta\label{5.0weak}
\end{equation} 
 as $h\rar 0$,  where $c=\( 64 / 3 \)^{ 1/2}$. 
 
 Equivalently
\begin{equation} { \int ( L^{ x+1}_{t}- L^{ x}_{ t})^{ 2}\,dx- 4t\over t^{ 3/4}}
\stackrel{\mathcal{L}}{\rightarrow}c\sqrt{\al_{1}}\,\,\eta\label{5.0tweak}
\end{equation}
as $t\rar\ff$.
\end{theorem} 

The equivalence of (\ref{5.0weak}) and (\ref{5.0tweak}) 
 follows from the scaling   relationship
 \begin{equation}
\{ L^{ x}_{ h^{-2}t}\,;\,( x,t)\in R^{ 1}\times R^{ 1}_{ +}\}
\stackrel{\mathcal{L}}{=}\{h^{ -1} L^{h x}_{ t}\,;\,( x,t)\in R^{ 1}\times
R^{ 1}_{ +}\},\label{scale}
\end{equation}
see e.g. \cite[Lemma 10.5.2]{book},
which implies that  
\begin{equation}
\int ( L^{ x+h}_{t}- L^{ x}_{ t})^{ 2}\,dx\stackrel{\mathcal{L}}{=}h^{3} \int ( L^{ x+1}_{t/h^{2}}- L^{ x}_{ t/h^{2}})^{ 2}\,dx.\label{scl}
\end{equation}
  Using this, and (\ref{5.0weak}) with $t=1$, and the change of variables $h^{2}=1/t$ gives (\ref{5.0tweak}).

\medskip	
  We show in Lemma \ref{lem-vare}  that   
\begin{equation}
E\(\int ( L^{ x+1}_{ t}- L^{ x}_{ t})^{ 2}\,dx\)=4  t +O( t^{1/2}).\label{9.13}
\end{equation}
 Consequently, (\ref{5.0tweak}) can be written as
\begin{equation} \qquad{ \int ( L^{ x+1}_{t}- L^{ x}_{ t})^{ 2}\,dx- E\(\int ( L^{ x+1}_{ t}- L^{ x}_{ t})^{ 2}\,dx\)\over t^{ 3/4}}
\stackrel{\mathcal{L}}{\rightarrow}c\sqrt{\al_{1}}\,\,\eta\label{5.0tweake}
\end{equation}
A similar weak law holds for (\ref{5.0weak}).

\medskip	
 In the proof of Theorem  \ref{theo-clt2}  we use the following result which is of independent interest: 
Let  $\{L^{ x }_{ t},\,\wt L^{ x }_{ t}\,;\,(x,t)\in R^{ 1}\times  R^{  1}_{ +}\}$ denote the local times of two independent  Brownian motions and let 
\begin{equation}
\bb_{s,t}= \int  L^{ x}_{ s}  \wt L^{ x}_{ t} \,dx\label{f7.1}
\end{equation}
denote their intersection local time.

\begin{theorem}\label{theo-cltindf} For each fixed $s,t$  
\begin{equation} { \int ( L^{ x+h}_{s}- L^{ x}_{ s})( \wt L^{ x+h}_{t}- \wt L^{ x}_{ t}) \,dx \over h^{ 3/2}}
\stackrel{\mathcal{L}}{\rightarrow}\wt C\sqrt{\bb_{s,t}}\,\,\eta\label{f5.0weaki}
\end{equation}  
as $h\rar 0$, where $\wt C=\({32 /3}\)^{ 1/2}$. Consequently
\begin{equation} { \int ( L^{ x+1}_{t}- L^{ x}_{ t})(\wt  L^{ x+1}_{t}- \wt L^{ x}_{ t}) \,dx \over t^{ 3/4}}
\stackrel{\mathcal{L}}{\rightarrow}\wt C\sqrt{\bb_{1,1}}\,\,\eta\label{f5.0tweaki}
\end{equation}
as $t\rar\ff$.
\end{theorem}

\medskip	
  We were motivated to try to find a central limit theorem for $\int   (L^{ x+h}_{ t} -L^{ x }_{
t})^{ 2}\,dx$ by our interest in
   the expression
\begin{equation} H_{ n}=\sum_{ i,j=1,\,i\neq j}^{ n}1_{ \{S_{ i}=S_{  j} \}}-{1
\over 2}\sum_{ i,j=1,\,i\neq j}^{ n}1_{
\{|S_{ i}-S_{ j}|=1 \}},\label{rp5c.4}
\end{equation} which appears as the Hamiltonian in a model for a polymer in a
repulsive medium,
\cite{HK}.  Here   $S:=\{S_{ n} \,;\,n=0,1,2,\ldots\}$ is a simple random walk on
$Z^{ 1}$. Note that
\be
H_{ n}=\frac{1}{2}\sum_{x\in Z^{ 1}}\(l_{ n}^{ x}-l_{ n}^{ x+1}\)^{ 2}
\ee
 where $l_{ n}^{
x}=\sum_{ i=1}^{ n}1_{ \{S_{ i}=x \}}$
 is the local time for  
$S$.  

\medskip	 Theorems \ref{theo-clt2} and \ref{theo-cltindf} are proved by the method of moments. In Section \ref{sec-ind} we show that Theorem \ref{theo-cltindf}   follows  immediately from moment estimates in Lemma \ref{lem-2.1}. Lemma \ref{lem-2.1} itself follows from  Lemma  \ref{lem-expind},   which obtains  the moments of an  expression analogous to the one in Lemma \ref{lem-2.1}, except that the fixed time $t$ is replaced by independent exponential times.  Lemma \ref{lem-expind} is proved in Section \ref{sec-expind}. Lemma \ref{lem-2.1} also requires Lemma \ref{lem-Lap} which allows us to use Laplace transform methods. Lemma \ref{lem-Lap} in proved in Section \ref{sec-Lap}. In Section \ref{sec-est} we derive some estimates on the potential densities of Brownian motion that  are used throughout this paper. In Section \ref{sec-BM} we show that Theorem \ref{theo-clt2}   follows from  Lemma \ref{lem-6.2}, on the moments of an  expression analogous to the left hand side of 
 (\ref{5.0weak}),  in which   $t$ is replaced by an independent exponential time.  
 Lemma \ref{lem-6.2} is proved in Section \ref{sec-clt}.   In Section \ref{sec-8} we obtain (\ref{9.13}).
 
  \medskip	 The basic tool we use for studying moments of local times is  Kac's moment formula. We use exponential times to make  Kac's moment formula manageable. Moments at exponential times correspond to the Laplace  transforms of the moments at fixed times. Since the left hand side of 
 (\ref{5.0weak}) has no obvious  monotonicity properties,   an important part of our proof involves showing how to derive  limit results for the moments of  (\ref{5.0weak}) from  limit results for their
 Laplace transforms.
 
 An alternate approach to proving Theorems \ref{theo-clt2} and \ref{theo-cltindf} would be to use  Tanaka's formula and martingale methods;  (see  \cite{Yor,YW}).     For the results in this paper   this would involve establishing   results about the differentiability of triple intersection local  times, as is done in  \cite{dsilt} for ordinary  intersection local  times. We   plan to return to this at a later date.
 
 \medskip	\noindent {\bf Acknowledgment:  }We thank Andrew Poje for numerical integrations which convinced us that  the results in Theorem \ref{theo-clt2} were correct and encouraged us to find a proof.

\section{Proof of Theorem \ref{theo-cltindf}}\label{sec-ind}

  We derive Theorem \ref{theo-cltindf} from the next lemma which is proved in this section.

\bl \label{lem-2.1}For all $s,t\geq 0$ and all integers $m\ge 0$
\begin{eqnarray} &&
\lim_{ h\rar 0}E\(\({  \int ( L^{ x+h}_{s}- L^{ x}_{ s})( \wt L^{ x+h}_{t}- \wt L^{ x}_{ t}) \,dx\over h^{ 3/2}}\)^{m}\)\nn\\ 
&&\hspace{ .5in}  =\left\{\begin{array}{ll}
\displaystyle{( 2n)!\over 2^{ n}n!}\( \displaystyle{32  \over 3}\)^{ n} E\lc\(\int L^{ x}_{ s}\wt L^{ x}_{ t} \,dx\)^{
n}\rc &\mbox{ if }m=2n\\\\
0&\mbox{ otherwise.}
\end{array}
\right.
\label{57.53}
\end{eqnarray}
\el

\noindent {\bf  Proof of Theorem \ref{theo-cltindf}  }\label{page 5}
 It follows from \cite[(6.12)]{CLR}   that   
 \begin{equation}
 E\lc\(\int (L^{ x}_{ s})^{ 2}\,dx\)^{ n}\rc\leq C_{s}^{ n}( (2n)!)^{ 1/4}.\label{57rb.1}
 \end{equation}
Therefore,
the right-hand side of (\ref{57.53}), which is the $2n$--th moment of $c\sqrt{\bb_{s,t}}\,\,\eta$ is less than or equal to $  C_{s,t}^{ n}( (2n)!)^{ 3/4}$. 
This implies that  $c\sqrt{\bb_{s,t}}\,\,\eta$ is determined by its moments; (see \cite[p. 227-228]{Feller}). Thus (\ref{f5.0weaki})  
follows  from  \cite[Theorem 30.2]{B},   which is often referred to as the method of moments.   We then get  (\ref{f5.0tweaki})   by using the  scaling relationship, (\ref{scale}).\qed

The next two lemmas are used in the proof of Lemma  \ref{lem-2.1}. Lemma \ref{lem-expind} is proved in Section \ref{sec-expind} and Lemma \ref{lem-Lap} is proved in Section \ref{sec-Lap}.

\medskip	Let $\la_{\ze}$ and $\wt\la_{\ze'}$ be independent  exponential times with  means $1/\ze$ and $1/\ze'$ respectively.

 \begin{lemma}\label{lem-expind}
For each integer $m\ge 0$, and   any $\ze,\ze'>0$,
\begin{eqnarray} &&
\lim_{ h\rar 0}E\(\( {\int ( L^{ x+h}_{\la_{\ze}}- L^{ x}_{ \la_{\ze}})(\wt L^{ x+h}_{\wt\la_{ \ze'}}- \wt L^{ x}_{\wt \la_{\ze'}}) \,dx \over h^{ 3/2} }  \)^{m}\)= a_{m}
\label{57.1}
\end{eqnarray}
where
\begin{equation}
a_{m} =\left\{\begin{array}{ll}
\displaystyle {( 2n)!\over 2^{ n}n!}\( {32 \over 3}\)^{ n} E\lc\(\int  L^{ x}_{ \la_{\ze}}  \wt L^{ x}_{\wt \la_{\ze'}}    \,dx\)^{
n}\rc &\mbox{ if }m=2n\\\\
0&\mbox{ otherwise.}
\end{array}
\right.\label{57.1a}
\end{equation}
\end{lemma}

We   write the statement of Lemma  \ref{lem-expind} in the form
\begin{eqnarray}&&
\lim_{ h\rar 0}\int_{0}^{\ff}\int_{0}^{\ff}e^{- \ze s- \ze' t} E\(\({ \int ( L^{ x+h}_{s}- L^{ x}_{s})(\wt L^{ x+h}_{t}- \wt L^{ x}_{t})\,dx \over h^{ 3/2}} \)^{m}\) \,ds\,dt
\nn\\
&&\qquad= \int_{0}^{\ff}\int_{0}^{\ff}e^{- \ze s- \ze' t}E\lc\eta^{m}\(  {32 \over 3}  \int  L^{ x}_{s}  \wt L^{ x}_{t}    \,dx \,\)^{
m/2}\rc \,ds\,dt\label{57.1d} .
\end{eqnarray}
 For   $h>0$ let 
\begin{eqnarray} &&
F_{h}(s,t;m) := E\(\({ \int ( L^{ x+h}_{s}- L^{ x}_{s})(\wt L^{ x+h}_{t}- \wt L^{ x}_{t})\,dx \over h^{ 3/2}} \)^{m}\) \label{57.1e},
\end{eqnarray}
and
\begin{eqnarray} &&
F_{0}(s,t;m) = E\lc\eta^{m}\(  {32 \over 3}  \int  L^{ x}_{s}  \wt L^{ x}_{t}    \,dx \,\)^{
m/2}\rc. \label{57.1e}
\end{eqnarray}
In this notation Lemma \ref{lem-expind} states that for  any $\ze,\ze'>0$ 
\begin{equation}
\lim_{ h\rar 0}\int_{0}^{\ff}\int_{0}^{\ff}e^{- \ze s- \ze' t} F_{h}(s,t;m) \,ds\,dt=
\int_{0}^{\ff}\int_{0}^{\ff}e^{- \ze s- \ze' t} F_{0}(s,t;m) \,ds\,dt.\label{57.1dd}
\end{equation}
(Note that $F_{0}(s,t;m) =0$ when $m$ is odd.)

\begin{lemma}\label{lem-Lap}  For  all  integers $m\ge 0$, and $0\leq h\leq 1$, $F_{h}(s,t;m)$ is a non--negative, polynomially bounded, continuous  increasing function of $(s,t)$.
\end{lemma}

\noindent{\bf Proof of Lemma  \ref{lem-2.1}  } 	
It follows from Lemma \ref{lem-Lap} that $F_{h}(s,t;m)$ is the distribution function of a measure $\mu_{h,m}$ on $R^{2}_{+}$; i.e.
\begin{equation}
F_{h}(s,t;m)=\int_{0}^{s}\int_{0}^{t}d\mu_{h,m}(u,v). \label{57.2}
\end{equation}
For any $0\le h\le 1$, consider
\begin{equation}
  \int_{0}^{\ff}\int_{0}^{\ff}e^{- \ze s- \ze' t} F_{h}(s,t;m) \,ds\,dt. 
   \end{equation}
Since $F_{h}(0,t;m)=F_{h}(s,0;m)=0 $,  it follows from integrating by parts, (in which we use   Lemma \ref{lem-Lap}),  that for all $\ze,\ze'>0$,
\begin{equation}
\ze \ze' \int_{0}^{\ff}\int_{0}^{\ff}e^{- \ze s- \ze' t}F_{h}(s,t;m)\,ds\,dt=\int_{0}^{\ff}\int_{0}^{\ff}e^{- \ze s- \ze' t}d\mu_{h,m}(s,t).\label{57.3}
\end{equation} 
We see from  (\ref{57.1dd}) and (\ref{57.3}) that for any $\ze,\ze'>0$,
\begin{equation}
\lim_{ h\rar 0}\int_{0}^{\ff}\int_{0}^{\ff}e^{- \ze s- \ze' t} d\mu_{h,m}(s,t)=
\int_{0}^{\ff}\int_{0}^{\ff}e^{- \ze s- \ze' t} d\mu_{0,m}(s,t).\label{57.1ddd}
\end{equation}

 It then follows from (\ref{57.1ddd}) and  the extended continuity theorem, \cite[Theorem 5.22]{K} that $\mu_{h,m}\stackrel{w}{\rar}\mu_{0,m}$. Using this and  Lemma \ref{lem-Lap} we see  that
\begin{equation}
\lim_{ h\rar 0}F_{h}(s,t)=F_{0}(s,t),\hspace{.2 in}\forall s,t,\label{57.4}
\end{equation}
which gives  (\ref{57.53}). (Actually,  \cite[Theorem 5.22]{K} is stated for probability measures on $R^{d}_{+}$. The case of general measures  on $R^{d}_{+}$ can be derived as in the proof of 
\cite[XIII.1, Theorem 2a]{Feller}. Unfortunately \cite[XIII.1, Theorem 2a]{Feller} only considers  measures  on $R^{1}_{+}$. Its extension to $R^{d}_{+}$ is routine.)

\qed

\section{Estimates for  the   potential density of \newline Brownian motion}\label{sec-est}

The $\al$-potential density of Brownian motion,
 \begin{equation}
u^{\al}(x)=\int_{0}^{\ff}e^{-\al t}p_{t}(x)\,dt={e^{-\sqrt{2\al}|x|} \over \sqrt{2\al}}\label{pot.1w}.
\end{equation}
Let $\la_{\al}$ be an independent exponential random variable with mean $1/\al$.

 Kac's moment formula, \cite[Theorem 3.10.1]{book}, states that
\begin{equation} E^{ x_{ 0}}\(\prod_{ j=1}^{ n}L^{ x_{ j}}_{ \la_{\al}} \)=\sum_{
\pi}\prod_{ j=1}^{ n}u^{\al}( x_{\pi( j)}-x_{\pi( j-1)})\label{1.2w}
\end{equation} where the sum runs over all permutations $\pi$ of $\{ 1,\ldots,
n\}$ and  
$\pi(0)=0.$

Let $\De_{ x}^{ h}$  denote
the finite difference operator on the variable $x$, i.e.
\begin{equation}
\De_{x}^{ h}\,f(x)=f(x+h)-f(x).\label{pot.3w}
\end{equation}
We write $\De^{ h}$ for $\De_{x}^{ h}$ when the variable $x$ is clear.

\medskip	The next lemma collects some facts about $u^{\al}(x)$ that are used in this paper.

\begin{lemma}\label{lem-vprop}Fix $\al,\bb>0$. For $0<h\leq 1$,
\bea
\De_{ x}^{ h}\De_{ y}^{ h} u^{\al}(x-y)\Bigg\vert_{ y=x}&= &2\({1-e^{-\sqrt{2\al}\,h} \over \sqrt{2\al}}\)=2h+O( h^{ 2}),\quad\label{1.8}
\\\nn\\
|\De ^{ h}\,u^{\al}(x)|&\leq &Ch\, u^{\al}( x),\label{1.3x}
\\\nn\\
|\De^{ h}\De^{ -h} u^{\al}(x )|&\leq& C h \, u^{\al}( x ),  \label{1.3yq}
\\\nn\\
|\De^{ h}\De^{ -h} u^{\al}(x )|&\leq& C h^{2}\, u^{\al}( x ), \hspace{.2 in}\forall \,|  x|\geq h.\label{1.3y}
\eea
In addition
\bea
\int  \(\De^{ h}\De^{ -h}\,u^{\al}(x)\) \(\De^{ h}\De^{ -h}\,u^{\bb}(x)\) \,dx&=&( 8/3+O( h))h^{ 3},\label{1.30g}
\\\nn\\ 
\int_{|x|\geq h} \(\De^{ h}\De^{ -h}\,u^{\al}(x)\)^{2}\,dx&=& O( h^{ 4}),\label{1.30gb}
\\\nn\\
\int |\De^{ h}\De^{ -h}\,u^{\al}(x)| \,dx&= &O( h^{ 2}).\label{li.13}
\eea
In all these statements the constants $C$ and the terms $O( h^{ \cd})$ may depend on $\al$ and $\bb$.
\end{lemma} 

\Proof
Since   
\begin{eqnarray} \lefteqn{
\De_{ x}^{ h}\De_{ y}^{ h} u^{\al}(x-y)\label{1.8w}}\\ && =
\{u^{\al}(x-y)-u^{\al}(x-y-h)\}\nonumber-
\{u^{\al}(x-y+h)-u^{\al}(x-y)\}\nonumber
\end{eqnarray}  
we have \begin{eqnarray} &&
\De_{ x}^{ h}\De_{ y}^{ h} u^{\al}(x-y)\Bigg\vert_{ y=x} =
\{u^{\al}(0)-u^{\al}(-h)\}-
\{u^{\al}(h)-u^{\al}(0)\}\nn\\ &&\hspace{ 1in}=2(u^{\al}(0)-u^{\al}(h) )=
2\({1-e^{-\sqrt{2\al}\,h} \over \sqrt{2\al}}\),\label{1.8a}
\end{eqnarray}
which gives (\ref{1.8}).

  To obtain  (\ref{1.3x}) we note that
\begin{equation}
\De_{x}^{ h}\,u^{\al}(x)=\({e^{-\sqrt{2\al}|x+h|} -e^{-\sqrt{2\al}|x|} \over \sqrt{2\al}}\).\label{pot.3ow}
\end{equation}
Therefore 
\bea
|\De_{x}^{ h}\,u^{\al}(x)|&=&{ e^{-\sqrt{2\al}|x |}\over \sqrt{2\al}}\left| e^{ \sqrt{2\al}(|x|-|x+h|)}-1\right| \label{pot.3owa}\\
&\le&  e^{-\sqrt{2\al}|x |}\( ||x|-|x+h||+O( ||x|-|x+h||^{2})\) \nn,
\eea
which gives  (\ref{1.3x}), (since we allow $C$ to depend on $\al$).
 
  To obtain  (\ref{1.3yq}) we note that  
 \bea
 |\De^{ h}\De^{- h}\,u^{\al}(x)|&=& |2u^{\al}(x) -  \,u^{\al}(x+h) -  \,u^{\al}(x-h) |\\
&\leq & |  \De^{ h}\,u^{\al}(x) | +|  \De^{ h}\,u^{\al}(x-h) |\nn,
 \eea
 and use (\ref{1.3x}).

To obtain  (\ref{1.3y}) we simply note that
when $|x|\geq h$, 
\begin{eqnarray} 
\De^{ h}\De^{ -h}\,u^{\al}( x)&=& 2u^{\al}( x)-u^{\al}( x+h)-u^{\al}( x-h)\label{1.26gd}\\ 
&=& u^{\al}( x)\(2-e^{-\sqrt{2\al}\,h}-e^{\sqrt{2\al}\,h}\)\nn .
\end{eqnarray}
The statement in  (\ref{1.30gb}) follows trivially from (\ref{1.3y}).

For (\ref{1.30g}) we note that for $|x|\leq h$
\begin{eqnarray} 
&&
\De^{ h}\De^{ -h}\,u^{\al}( x)\label{1.26g}\\&&= 2u^{\al}( x)-u^{\al}( x+h)-u^{\al}( x-h)\nn\\ &&= (u^{\al}( 0)-u^{\al}(
x+h))+(u^{\al}( 0)-u^{\al}( x-h))-2(u^{\al}( 0)-u^{\al}( x))\nn\\ &&= | x+h|+ | x-h|-2 | x|+O( h^{ 2}).\nn
\end{eqnarray} 
 Therefore when $0\leq x\leq h$, we  have 
\begin{equation}
\De^{ h}\De^{ -h}\,u^{\al}( x)=x+h+h-x-2x+O( h^{ 2})=( 2+O( h))(h-x)\label{1.27g}
\end{equation} 
  and similarly for $\De^{ h}\De^{ -h}\,u^{\bb}( x)$. Consequently
\bea
\int_{0}^{ h} \(\De^{ h}\De^{ -h}\,u^{\al}(x)\) \(\De^{ h}\De^{ -h}\,u^{\bb}(x)\) \,dx \nn
&=& ( 4+O( h))\int_{0}^{
h}(h-x)^{2}\,dx\\
&=&( 4/3+O( h))h^{ 3}.\label{1.28g}
\eea
 Similarly, when $-h\leq x\leq 0$ it follows from  (\ref{1.26g}) that
\begin{equation}
\De^{ h}\De^{ -h}\,u^{\al}( x)=h-x+x+h+2x+O( h^{ 2})=( 2+O( h))(h+x)\label{1.29g}
\end{equation} 
and similarly for $\De^{ h}\De^{ -h}\,u^{\bb}( x)$. Consequently
\bea \int_{-h}^{ 0}\(\De^{ h}\De^{ -h}\,u^{\al}(x)\)\(\De^{ h}\De^{ -h}\,u^{\bb}(x)\) \,dx&=& ( 4+O( h))\int_{-h}^{
0}(h+x)^{2}\,dx\nn\\
&=&( 4/3+O( h))h^{ 3}.\label{1.30gx}
\eea
Using  (\ref{1.28g}), (\ref{1.30gx}) and   (\ref{1.30gb}) we get (\ref{1.30g}).

To obtain (\ref{li.13}) we write 
\begin{eqnarray}
&&\int |\De^{ h}\De^{- h}\,u^{\al}(y) |\,dy 
\label{li.13a}\\
&&\qquad= \int_{|y|\leq h} |\De^{ h}\De^{- h}\,u^{\al}(y) |\,dy   + \int_{|y|\geq h} |\De^{ h}\De^{- h}\,u^{\al}(y) |\,dy \nonumber\\
&&\qquad\leq Ch\int_{|y|\leq h} 1\,dy   + Ch^{2}\int_{|y|\geq h}  u^{\al}(y) \,dy=O( h^{
2}) ,\nonumber
\end{eqnarray}
where for the last line we use (\ref{1.3yq}) and (\ref{1.3y}).  \qed

\section{Proof of Lemma \ref{lem-expind} }\label{sec-expind}

Let $X_{t},\, \wt X_{t}$ be two independent  Brownian motions in $R^{1}$. Let  $L_{t}^{x},\, \wt L^{x}_{t}$ 
denote their local times, and  let $\la_{\ze},\wt\la_{\ze'}$ be independent  exponential times of mean $1/\ze,1/\ze'$ respectively. Set
\begin{equation}
\bb_{2}=\int L_{\la_{\ze}}^{x}\,\wt L^{x}_{\wt\la_{\ze'}}\,dx.\label{b.1}
\end{equation}

It follows from (\ref{1.2w}),  the Kac  moment formula,  that 
\begin{eqnarray}
&& 
E\(\prod_{i=1}^{m}      L^{ x_{i}}_{\la_{\ze}}\, \wt  L^{ y_{i}}_{\wt\la_{\ze'}}  \) = 
E\(\prod_{i=1}^{m}      L^{ x_{i}}_{\la_{\ze}}   \) E\(\prod_{i=1}^{m}  \wt  L^{ y_{i}}_{\wt \ze}  \)\label{f7.2kac}\\
&&\qquad= \sum_{\pi}  \prod_{j=1}^{m}u^{\ze}(x_{\pi(j)}-x_{\pi(j-1)})\,dr_{j}\nonumber\\
&&\qquad\qquad \times \sum_{\pi'}  \prod_{j=1}^{m}u^{\ze'} (y_{\pi'(j)}-y_{\pi'(j-1)}),
\nonumber
\end{eqnarray}
where the sums run over all permutations  $\pi$ and $\pi'$ of $\{1,\ldots, m\}$,    $\pi(0) =\pi'(0) =0$
and $x_{0} =0$.
Consequently, by setting each $y_{i}$ equal to $x_{i}$ we see that
\begin{eqnarray}
&&E\(   \(\int L^{ x }_{  \la_{\ze}}\, \wt  L^{x}_{\wt\la_{\ze'}}\,dx \)^{m}\)=E\(  \prod_{i=1}^{m} \int L^{ x_{i}}_{ \la_{\ze}}\, \wt  L^{ x_{i}}_{\wt\la_{\ze'}}\,dx_{i}  \)
\label{kac.m}\\
&& \qquad= \sum_{\pi,\pi'} \int\( \prod_{j=1}^{m}u^{\ze}(x_{\pi(j)}-x_{\pi(j-1)})\right. \nonumber\\
&&\left.\qquad \qquad \times   \prod_{j=1}^{m}u^{\ze'}(x_{\pi'(j)}-x_{\pi'(j-1)})\, \) \prod_{i=1}^{m}\,dx_{i}.
\nonumber
\end{eqnarray}
Similarly
\begin{eqnarray}
&&
E\(\prod_{i=1}^{m}  ( L^{ x_{i}+h}_{\la_{\ze}}- L^{ x_{i}}_{  \la_{\ze}})(\wt   L^{ y_{i}+h}_{\wt\la_{\ze'}}- \wt  L^{ y_{i}}_{\wt\la_{\ze'}})  \)  \label{f7.2}\\
&&\qquad=\(\prod_{i=1}^{m}\Delta_{x_{i}}^{h}\Delta_{y_{i}}^{h}\)   
E\(\prod_{i=1}^{m}      L^{ x_{i}}_{\la_{\ze}}\, \wt  L^{ y_{i}}_{\wt\la_{\ze'}}  \) \nonumber\\
&&\qquad=\(\prod_{i=1}^{m}\Delta_{x_{i}}^{h}\Delta_{y_{i}}^{h}\)   
E\(\prod_{i=1}^{m}      L^{ x_{i}}_{ \la_{\ze}}   \) E\(\prod_{i=1}^{m}  \wt  L^{ y_{i}}_{\wt\la_{\ze'}}  \)\nonumber\\
&&\qquad=\(\prod_{i=1}^{m}\Delta_{x_{i}}^{h}\) \sum_{\pi} \prod_{j=1}^{m}u^{\ze}(x_{\pi(j)}-x_{\pi(j-1)}) \nonumber\\
&&\qquad\qquad \times \(\prod_{i=1}^{m}\Delta_{y_{i}}^{h}\) \sum_{\pi'}  \prod_{j=1}^{m}u^{\ze'}(y_{\pi'(j)}-y_{\pi'(j-1)}) .
\nonumber
\end{eqnarray} 

 Using the product rule for finite differences, 
\begin{equation}
\De^{ h}(fg)(x)=(\De^{ h}f(x))g(x+h)+f(x)\De^{ h}g(x)\label{prf}
\end{equation}
 we can write 
\bea 
\lefteqn{ \(\prod_{i=1}^{m}\Delta_{x_{i}}^{h}\) \sum_{\pi} \prod_{j=1}^{m}u^{\ze}(x_{\pi(j)}-x_{\pi(j-1)}) \nn}\\
  &&\, =\sum_{ \pi, a }  \prod_{ j=1}^{ m}\(\(\De^{ h}_{ x_{ \pi( j)}}\)^{a_{ 1}(j)}
\(\De^{ h}_{ x_{ \pi( j-1)}}\)^{a_{ 2}(j)}\,u^{\ze,\sharp}(x_{\pi(j)}-x_{\pi(j-1)})\)  \label{3.2}
 \eea
where   the sum runs over   $\pi $  and   all   
$a =(a_{ 1},a_{ 2})\,:\,[1,\ldots, m]\mapsto \{ 0,1\} \times \{ 0,1\} $, with the
restriction that for each $i$ there is exactly one  factor  of the form $\De^{
h}_{ x_{i}}$.  (Here we define $(\De_{x_{i}}^{h})^{0}=1 $ and $(\De_{0}^{h}) =1 $.) In this formula, $u^{\ze,\sharp}(x)$ can take any of the values $u^{\ze}(x)$, $u^{\ze}(x+h)$ or $u^{\ze}(x-h)$. (We consider all three possibilities in the subsequent proofs.) It is important to recognize that in (\ref{3.2}) each of  the difference operators is applied to only one of the  terms   $u^{\ze,\sharp}(\cdot)$.

Using (\ref{3.2}) we see that if we set $x_{i}=y_{i}$, $i=0,\ldots,m$ in (\ref{f7.2})  we get 
\bea && E\(\( \int ( L^{ x+h}_{\la_{\ze}}- L^{ x}_ {\la_{\ze}})(\wt   L^{ x+h}_{ \wt\la_{\ze'}}-\wt   L^{ x}_{\wt\la_{\ze'}}) \,dx \)^{ m}
\)\label{f1.20gi}\\ &&\qquad= \sum_{ \pi,\pi',a,a'}\int \mathcal{T}'_{h}( x;\,\pi,\pi',a,a')\,dx\nn
\eea 
where   $x=(x_{1},\ldots, x_{m})$ and 
\begin{eqnarray}
\lefteqn{ 
\mathcal{T}'_{h}( x;\,\pi,\pi',a,a')\label{3.9}}\\
&&\qquad=   \prod_{ j=1}^{ m}\(\(\De^{ h}_{ x_{ \pi( j)}}\)^{a_{ 1}(j)}
\(\De^{ h}_{ x_{ \pi( j-1)}}\)^{a_{ 2}(j)}\,u^{\ze,\sharp}(x_{\pi(j)}-x_{\pi(j-1)})\)\, \nn \\
&&\hspace{.5 in}  \prod_{ k=1}^{ m}\(\(\De^{ h}_{ x_{ \pi'( k)}}\)^{a'_{ 1}(k)}
\(\De^{ h}_{ x_{ \pi'( k-1)}}\)^{a'_{ 2}(k)}\,u^{\ze',\sharp}( x_{\pi'(k)}- x_{\pi'(k-1)})\),\nn
\end{eqnarray}
and  where the sum runs over all permutations  $\pi$  and $\pi'$ of $\{1,\ldots, m\}$ and   all  
$a=(a_{ 1},a_{ 2})\,:\,[1,\ldots, m]\mapsto \{ 0,1\}\times \{ 0,1\}$ and $a'=(a'_{ 1},a'_{ 2})\,:\,[1,\ldots, m]\mapsto \{ 0,1\}\times \{ 0,1\}$ with the
restriction that for each $i$ there   is   exactly one  factor  of the form $\De^{
h}_{ x_{i}}$ in the second line of   (\ref{3.9}), and similarly  in the third line of (\ref{3.9}).  

Let
\begin{eqnarray}
\lefteqn{ 
\mathcal{T}_{h}( x;\,\pi,\pi',a,a') \nn}\\
&&\qquad= \prod_{ j=1}^{ m}\(\(\De^{ h}_{ x_{ \pi( j)}}\)^{a_{ 1}(j)}
\(\De^{ h}_{ x_{ \pi( j-1)}}\)^{a_{ 2}(j)}\,u^{\ze} (x_{\pi(j)}-x_{\pi(j-1)})\)\label{f1.21gi} \\
&&\hspace{.15 in}\prod_{ k=1}^{ m}\(\(\De^{ h}_{ x_{ \pi'( k)}}\)^{a'_{ 1}(k)}
\(\De^{ h}_{ x_{ \pi'( k-1)}}\)^{a'_{ 2}(k)}\,u^{\ze'} ( x_{\pi'(k)}- x_{\pi'(k-1)})\).\nn
\end{eqnarray}
The difference between (\ref{3.9}) and (\ref{f1.21gi}) is that   $u^{\ze,\sharp}$
is replaced by $u^{\ze}$ and similarly for $u^{\ze',\sharp}$. To simplify the computations we first obtain   
\begin{equation}
   \lim_{h\to 0}\frac{1}{h^{3m/2}} \sum_{ \pi,\pi',a,a'}\int \mathcal{T}_{h}( x;\,\pi,\pi',a,a')\,dx\label{3.11}
   \end{equation} 
   and then explain why (\ref{3.11}) is unchanged when $\TT_{h}$ is replaced by $\TT'_{h}$
 
\medskip	
 Recall that $\De^{h}f(u)=f(u+h)-f(h)$ so that 
 \begin{equation}
   \De^{h}\De^{-h}f(u-v)=2f(u-v)-f(u-v-h)-f(u-v+h).
   \end{equation}
Consequently
 \begin{equation}
    \De_{u}^{h}\De_{v}^{h}f(u-v)=  \De^{h}\De^{-h}f(u-v) .\label{3.7}
      \end{equation}
      
\medskip	      We proceed to evaluate  (\ref{3.11}) based on the different ways the difference operators in (\ref{f1.21gi}) are distributed. We examine these in three subsections. The reader will see that   the only limits in (\ref{3.11}), that are not zero,  come from the case considered in  Subsection \ref{ss3.1}.

\medskip	
Let  $e=(e(1),\ldots,e(2n))$ where 
 $e(2j)=(1,1),\,e(2j-1)=(0,0)$,   $j=1\ldots n$.

\subsection{   ${\bf a=a'=e}$ and compatible permutations}\label{ss3.1}

Let $m=2n$  and let $\mathcal{P}=\{(l_{2i-1},l_{2i})\,,\,1\leq i\leq n\}$ be a pairing of the integers $[1,2n]$. Let $\pi$ and $\pi'$ be two permutations  of   $[1,2n]$   such that for each $1\leq j\leq n$, 
$\{\pi(2j-1), \pi(2j)\}=\{l_{2i-1},l_{2i}\}$ for some, necessarily unique, $ 1\leq i\leq n$ and similarly for $\pi'$,   i.e.  for each $1\leq j\leq n$,  
$\{\pi'(2j-1), \pi'(2j)\}=\{l_{2k-1},l_{2k}\}$ for some, necessarily unique, $ 1\leq k\leq n$. In this case we say that $\pi$ and $\pi'$ are compatible with the pairing $\mathcal{P}$ and write this  as $(\pi,\pi')\sim \mathcal{P}$. (Note that $\{\pi(2j-1), \pi(2j)\} $  is  not necessarily equal to  $\{\pi'(2j-1), \pi'(2j)\}$.    Furthermore, when we write $\{\pi(2j-1), \pi(2j)\}=\{l_{2i-1},l_{2i}\}$ we mean as two  sets, so, according to what $\pi$ is, we may have  $\pi(2j-1)=l_{2i-1}$, $\pi(2j )=l_{2i }$ or $\pi(2j-1)=l_{2i}$, $\pi(2j )=l_{2i-1 }$.   A similar situation exist for $\pi'$.) We write $\pi\sim \pi'$ to mean that 
$(\pi,\pi')\sim \mathcal{P}$ for some pairing $\mathcal{P}$. In this subsection we show that  
 \begin{eqnarray} \quad
&&\sum_{\pi\sim \pi'}\int \mathcal{T}_{h}( x;\,\pi,\pi',e,e) \,\prod_{ j=1}^{ 2n}\,dx_{j}\label{f9.45sj}\\
&& \qquad  =  {( 2n)!\over 2^{ n}n!} \({32h^{3}\over 3} \)^{n}  E\lc\(\int  L^{ x}_{ \la_{\ze}}  \wt L^{ x}_{\wt \la_{\ze'}}    \,dx\)^{
n}\rc+O(h^{3n+1}).\nn
 \end{eqnarray} 
  In   Subsections \ref{ss-3.2} and \ref{ss3.4} we show that
\begin{equation}
\sum_{\pi\not\sim \pi'\mbox{ or }(a,a')\neq (e,e)  }\Big |\int \mathcal{T}_{h}( x;\,\pi,\pi',a,a') \,\prod_{ j=1}^{ 2n}\,dx_{j}\Big |=O(h^{3n+1}). \label{big.2}
\end{equation}
 Together these   estimates give (\ref{57.53}).

 When $\pi$ and $\pi'$ are compatible   it follows from   (\ref{f1.21gi}) and (\ref{3.7}) that
\begin{eqnarray}
&& 
\mathcal{T}_{h}( x;\,\pi,\pi',e,e) \label{f1.21gia}\\
&&\qquad=  \prod_{ j=1}^{ n}\(\De^{ h}\De^{ -h}
\,u^{\ze}(x_{\pi(2j)}-x_{\pi(2j-1)})\) \nn\\ &&\hspace{
1.1in}\times\prod_{ j=1}^{ n}\, u^{\ze}(x_{\pi(2j-1)}-x_{\pi(2j-2)}) \nn\\
&&\qquad\hspace{.2 in}  \prod_{ k=1}^{ n}\(\De^{ h}\De^{ -h}
\,u^{\ze'}(x_{\pi'(2k)}-x_{\pi'(2k-1)})\) \nn\\ &&\hspace{
1.1in}\times\prod_{ k=1}^{ n}\, u^{\ze'}(x_{\pi'(2k-1)}-x_{\pi'(2k-2)}).\nn
\end{eqnarray}

We would like to integrate ${T}_{h}( x;\,\pi,\pi',e,e) $  with respect to $x$ but this is not easy because the variables 
\[\{x_{\pi(2j)}-x_{\pi(2j-1)}\,,x_{\pi'(2j)}-x_{\pi'(2j-1)}\,,\,j\in [1,n]\}\] and \[\{x_{\pi(2j-1)}-x_{\pi(2j-2)}\,,x_{\pi'(2j-1)}-x_{\pi'(2j-2)}\,,\,j\in [1,n]\}\]  are not independent.   To get around this difficulty we first write
\begin{equation}
1=\prod_{i=1}^{n}\(1_{\{|x_{l_{2i}}-x_{l_{2i-1}}|\leq h\}}
+1_{\{|x_{l_{2i}}-x_{l_{2i-1}}|\geq h\}}\)=\sum_{A\subseteq  [1,\ldots,n]}1_{D_{A}}\label{f9.31}
\end{equation}
where
\begin{equation}
D_{A}=  \{|x_{l_{2i}}-x_{l_{2i-1}}|\leq h,\,i\in A\} \cap  \{|x_{l_{2i}}-x_{l_{2i-1}}|> h,\,i\in A^{c}\}.\label{f9.31b}
\end{equation}
 and use it to write 
\begin{eqnarray}
\lefteqn{\int \mathcal{T}_{h}( x;\,\pi,\pi',e,e) \,\prod_{ j=1}^{ 2n}\,dx_{j}\label{f9.32}}\\
&&=\int \prod_{i=1}^{n}\(1_{\{|x_{l_{2i}}-x_{l_{2i-1}}|\leq \sqrt{h}\}}\)
\mathcal{T}_{h}( x;\,\pi,\pi',e,e) \prod_{ j=1}^{ 2n}\,dx_{j}+E_{1,h}\nn
\end{eqnarray}
where 
\begin{eqnarray}
&&E_{1,h}=\sum_{A^{ c}\ne \phi}\int_{ D_{A} } 
\mathcal{T}_{h}( x;\,\pi,\pi',e,e)\prod_{ j=1}^{ 2n}\,dx_{j}.
\label{f9.32a}
\end{eqnarray} 
 Let 
\be w^{\ze}(x)=|\De^{ h}\De^{ -h}
\,u^{\ze}(x)|.\label{4.20}
\ee We have 
\begin{eqnarray}
	\lefteqn{\left|\int_{ D_{A} } 
\mathcal{T}_{h}( x;\,\pi,\pi',e,e)\prod_{ j=1}^{ 2n}\,dx_{j}\right| 
\label{f9.34}}\\
&&\leq \int_{ D_{A} }  \prod_{ j=1}^{ n}\, u^{\ze}(x_{\pi(2j-1)}-x_{\pi(2j-2)})\, w^{\ze}(x_{\pi(2j)}-x_{\pi(2j-1)})\nn  \\  &&\hspace{
.4in}\times\prod_{ k=1}^{ n}\,u^{\ze'}(x_{\pi'(2k-1)}-x_{\pi'(2k-2)})\,w^{\ze'}(x_{\pi'(2k)}-x_{\pi'(2k-1)})   \, \prod_{ j=1}^{ 2n}\,dx_{j}\nn.
\eea

Let
\begin{equation}
\wt D_{A}=\{|x_{2j-1}-x_{2j-2}|\leq h,\,j\leq |A|\} \cap \{|x_{2j-1}-x_{2j-2}|> h,\,j> |A|\}\label{f9.31bv}
\end{equation}
Applying the Cauchy-Schwarz inequality in (\ref{f9.34}) to separate the   terms in $\pi$ from the terms in  $\pi'$, and then relabeling, we get
\bea
&&
\Bigg|\int_{ D_{A} } 
\mathcal{T}_{h}( x;\,\pi,\pi',e,e)\prod_{ j=1}^{ 2n}\,dx_{j}\Bigg|^{2}\nn\\
&&\qquad  \leq   \int_{\wt  D_{A} }   \prod_{ j=1}^{ n}\, (u^{\ze}(x_{2j-1}-x_{2j-2}))^{2}\,  \, (w^{\ze}(x_{2j}-x_{2j-1}))^{2}\, \prod_{ j=1}^{ 2n}\,dx_{j} \nn\\
&&\qquad\times  \int_{\wt  D_{A} }   \prod_{ j=1}^{ n}\, (u^{\ze'}(x_{2j-1}-x_{2j-2}))^{2}\,  \, (w^{\ze'}(x_{2j}-x_{2j-1}))^{2}\, \prod_{ j=1}^{ 2n}\,dx_{j} \nn\\
&&\qquad\leq  \(\prod_{j=1}^{| A|} \int    \,\, (w^{\ze}(x_{2j} ) )^{2}   \, \,dx_{2j}\)\(\prod_{j=| A|+1}^{n} \int     1_{\{|x_{2j}|\geq h\}}  \,\, (w^{\ze}(x_{2j} ) )^{2}   \, \,dx_{2j}\)\nn\\
&&\qquad\times  \(\prod_{j=1}^{| A|} \int    \,\, (w^{\ze'}(x_{2j} ) )^{2}   \, \,dx_{2j}\)\(\prod_{j=| A|+1}^{n} \int     1_{\{|x_{2j}|\geq h\}}  \,\, (w^{\ze'}(x_{2j} )  )^{2}  \, \,dx_{2j}\)\nonumber\\
&&\qquad\leq \(C^{n}h^{3n+|A^{c}|}\)^{2} ,\label{3.15}
\end{eqnarray} 
where 
the last inequality comes from (\ref{1.30g}) and (\ref{1.30gb}). Combining this with (\ref{f9.32a}) we see that  
\begin{equation}
|E_{1,h}|=O(h^{3n+1}).\label{3.15a}
\end{equation}

\medskip	 
We now study  
\begin{equation} \hspace{.4 in}
  \wt B_{h}(\pi,\pi',e,e ):=\int \prod_{i=1}^{n}\(1_{\{|x_{l_{2i}}-x_{l_{2i-1}}|\leq h\}}\)
\mathcal{T}_{h}( x;\,\pi,\pi',e,e)\prod_{ j=1}^{ 2n}\,dx_{j}. \label{f9.36}
\end{equation}
 Recall that for  each $1\le j\le n$, 
$\{\pi(2j-1), \pi(2j)\}=\{l_{2i-1},l_{2i}\} $, for some $1\le i\le n$.  We identify these relationships by setting $i=\si (j) $ when  $\{\pi(2j-1), \pi(2j)\}=\{l_{2i-1},l_{2i}\} $.     We write
\begin{eqnarray}
&&\prod_{ j=1}^{ n}\, u^{\ze}(x_{\pi(2j-1)}-x_{\pi(2j-2)})
\label{f9.37}\\
&&\qquad=\prod_{ j=1}^{ n}\,\( u^{\ze}(x_{l_{2\si (j)-1}}-x_{l_{2\si (j-1)-1}})+\De^{h_{j}}u^{\ze}(x_{l_{2\si (j)-1}}-x_{l_{2\si (j-1)-1}})\) ,   \nn
\end{eqnarray}
where $h_{j}=(x_{\pi(2j-1)}-x_{l_{2\si (j)-1}})+(x_{l_{2\si (j-1)-1}}-x_{\pi(2j-2)})$. 
Note that because of the presence of  the term $\prod_{i=1}^{n}\(1_{\{|x_{l_{2i}}-x_{l_{2i-1}}|\leq h\}}\)$ in the integral in (\ref{f9.36})  we need only be concerned with values of $|h_{j}|\leq 2h$, $1\le j\le n$.

  Similarly we set   $i=\si' (j)$ when
$\{\pi'(2j-1), \pi'(2j)\}=\{l_{2i-1},l_{2i}\} $, and write
\begin{eqnarray}
&&\quad\prod_{ j=1}^{ n}\, u^{\ze'}(x_{\pi'(2j-1)}-x_{\pi'(2j-2)})
\label{f9.37s}\\
&&\qquad=\prod_{ j=1}^{ n}\,\( u^{\ze'}(x_{l_{2\si' (j)-1}}-x_{l_{2\si' (j-1)-1}})+\De^{h'_{j}}u^{\ze'}(x_{l_{2\si' (j)-1}}-x_{l_{2\si' (j-1)-1}})\)   \nonumber
\end{eqnarray}
where $h'_{j}=(x_{\pi'(2j-1)}-x_{l_{2\si' (j)-1}})+(x_{l_{2\si' (j-1)-1}}-x_{\pi'(2j-2)})$.  As above  we need only be concerned with values of $|h'_{j}|\leq 2h$, $1\le j\le n$.
 
We substitute  (\ref{f9.37}) and  (\ref{f9.37s}) into the term $\mathcal{T}_{h}( x;\,\pi,\pi',e,e)$ in (\ref{f9.36}) and expand the products so that we can write 
 $\wt B_{h}(\pi,\pi',e,e )$ as a sum of many terms to get
\begin{eqnarray} 
\wt B_{h}(\pi,\pi',e,e )  &:=&\int \prod_{i=1}^{n}\(1_{\{|x_{l_{2i}}-x_{l_{2i-1}}|\leq h\}}\)
\wt \mathcal{T}_{h}( x;\,\pi,\pi',e,e )\prod_{ j=1}^{ 2n}\,dx_{j}+E_{2,h}\nn\\
&&\label{4.28}
\end{eqnarray}
where 
\begin{eqnarray}
\lefteqn{
\wt \mathcal{T}_{h}( x;\,\pi,\pi',e,e)  
 =  \prod_{ i=1}^{ n}\De^{ h}\De^{ -h}
\,u^{\ze}(x_{l_{2i}}-x_{l_{2i-1}})}\label{f1.21giaww}\\ &&\hspace{
2in}\times\prod_{ j=1}^{ n}\,  u^{\ze}(x_{l_{2\si (j)-1}}-x_{l_{2\si (j-1)-1}})\,\,  \nn\\
&&\hspace{1 in}  \prod_{ i=1}^{ n}\De^{ h}\De^{ -h}
\,u^{\ze'}(x_{l_{2i}}-x_{l_{2i-1}})\,  \nn\\ &&\hspace{
2in}\times\prod_{ j=1}^{ n}\,  u^{\ze'}(x_{l_{2\si' (j)-1}}-x_{l_{2\si' (j-1)-1}}) \nn
\end{eqnarray}
and
\begin{equation}
 E_{2,h}=\sum_{ {A, A'\subseteq  [1,\ldots,n]} } E_{2,h;A, A'}\label{f9.40r},
\end{equation}
where 
\begin{eqnarray}
\lefteqn{
E_{2,h;A, A'}:=\int \prod_{i=1}^{n}\(1_{\{|x_{l_{2i}}-x_{l_{2i-1}}|\leq h\}}\) \hspace{.2 in} \prod_{ j=1}^{ n}\De^{ h}\De^{ -h}
\,u^{\ze}(x_{\pi(2j)}-x_{\pi(2j-1)})\, \nn}\\
&&\hspace{
.1in}\times\prod_{ j\in A} \,  u^{\ze}(x_{l_{2\si (j)-1}}-x_{l_{2\si (j-1)-1}})\prod_{ j\in A^{c}} \, \De^{h_{j}} u^{\ze}(x_{l_{2\si (j)-1}}-x_{l_{2\si (j-1)-1}}) \nn\\
&&\hspace{1.3 in}\times\  \prod_{ k=1}^{ n}\De^{ h}\De^{ -h}
u^{\ze'}(x_{\pi'(2k)}-x_{\pi'(2k-1)})\label{f9.40}\\ && \times\prod_{ k\in A'} \,  u^{\ze'}(x_{l_{2\si' (k)-1}}-x_{l_{2\si' (k-1)-1}})\nn\\
&&\hspace{1.3in}\times\prod_{ k\in A'^{c}} \, \De^{h'_{k}} u^{\ze'}(x_{l_{2\si' (k)-1}}-x_{l_{2\si' (k-1)-1}}) \prod_{ j=1}^{ 2n}\,dx_{j},\nn
\end{eqnarray}
  and   $A^{c}$ and $A'^{c}$ are not both empty. 

Using (\ref{1.3x}) we see that 
\begin{eqnarray}
&& 
|E_{2,h;A, A'}|\leq C^{n} h^{|A^{c}|+|A'^{c}|}\int \prod_{i=1}^{n}\(1_{\{|x_{l_{2i}}-x_{l_{2i-1}}|\leq h\}}\) \label{f9.40a}\\
&&\hspace{.2 in} \prod_{ j=1}^{ n}w^{\ze}(x_{\pi(2j)}-x_{\pi(2j-1)})\,\prod_{ j=1}^{n} \,  u^{\ze}(x_{l_{2\si (j)-1}}-x_{l_{2\si (j-1)-1}}) \nn\\
&&\hspace{.2 in}  \prod_{ k=1}^{ n}w^{\ze'}(x_{\pi'(2k)}-x_{\pi'(2k-1)})\,\prod_{ k=1}^{n} \,  u^{\ze'}(x_{l_{2\si (k)-1}}-x_{l_{2\si (k-1)-1}})\,dx_{j}.\nn
\end{eqnarray}
Using the Cauchy-Schwarz inequality as   in (\ref{3.15}), we see that
\begin{equation}
|E_{2,h;A, A'}|\leq C h^{3n+|A^{c}|+|A'^{c}|}\label{f9.41z}.
\end{equation}
It   now follows from (\ref{f9.40r}) that
\begin{equation}
E_{2,h }=O(h^{3n+1}).\label{f9.41z2}
\end{equation}

\medskip	 
We now  consider 
\begin{equation}\qquad 
\ov B_{h}(\pi,\pi',e,e ):=\int \prod_{i=1}^{n}\(1_{\{|x_{l_{2i}}-x_{l_{2i-1}}|\leq h\}}\)
\wt \mathcal{T}_{h}( x;\,\pi,\pi',e,e )\prod_{ j=1}^{ 2n}\,dx_{j}. \label{9.42}
\end{equation}
  Using (\ref{f9.32})   and (\ref{3.15a}) we see that
\begin{equation} 
\ov B_{h}(\pi,\pi',e,e )=\int  
\wt \mathcal{T}_{h}( x;\,\pi,\pi',e,e )\prod_{ j=1}^{ 2n}\,dx_{j}+O(h^{3n+1}).\label{9.43}
\end{equation} 

  Using (\ref{f1.21giaww}) we see that 
\begin{eqnarray}
  \lefteqn{
\int  
\wt \mathcal{T}_{h}( x;\,\pi,\pi',e,e )\prod_{ j=1}^{ 2n}\,dx_{j}\nn}\\
&&\label{f9.44d}=\int  \prod_{ i=1}^{ n}\De^{ h}\De^{ -h}
\,u^{\ze}(x_{l_{2i}}-x_{l_{2i-1}})  \\ &&\hspace{
.7in}\times\prod_{ j=1}^{ n}\,  u^{\ze}(x_{l_{2\si (j)-1}}-x_{l_{2\si (j-1)-1}})\, \nn\\
&&\hspace{.3 in}  \prod_{ i=1}^{ n}\De^{ h}\De^{ -h}
\,  u^{\ze'}(x_{l_{2i}}-x_{l_{2i-1}})\ \nn\\ &&\hspace{
.7in}\times\prod_{ j=1}^{ n}\,   u^{\ze'}(x_{l_{2\si' (j)-1}}-x_{l_{2\si' (j-1)-1}}) \prod_{i=1}^{2n}\,dx_{i}\nn.
\eea
  We make the change of   variables $x_{l_{2i}}\to  x_{l_{2i}}+x_{l_{2i-1}}  $, $i=1,\ldots,n$ and  write this as  
\bea
&& \int    \prod_{ i=1}^{ n}\De^{ h}\De^{ -h}
\,u^{\ze}(x_{l_{2i}} )\,  \prod_{ j=1}^{ n}\,  u^{\ze}(x_{l_{2\si (j)-1}}-x_{l_{2\si (j-1)-1}})\, \nn\\
&&\hspace{.3 in}  \times  \prod_{ i=1}^{ n}\De^{ h}\De^{ -h}
\,  u^{\ze'}(x_{l_{2i}} ) \prod_{ j=1}^{ n}\,   u^{\ze'}(x_{l_{2\si' (j)-1}}-x_{l_{2\si' (j-1)-1}})  \prod_{i=1}^{2n}\,dx_{i}.\nn
\end{eqnarray}

We  now rearrange the integrals with respect to $x_{l_{2}} ,x_{l_{4}} ,\ldots,x_{l_{2n}} $ and get
\begin{eqnarray}
\lefteqn{
\int  
\wt \mathcal{T}_{h}( x;\,\pi,\pi',e,e )\prod_{ j=1}^{ 2n}\,dx_{j}=  \(\int   \(\De^{ h}\De^{ -h}
\,u^{\ze}(x )\) \(\De^{ h}\De^{ -h}
\,u^{\ze'}(x )\)  \,dx\)^{n}\nn}\\
&&   \times\int \prod_{ j=1}^{ n}\,    u^{\ze}(x_{l_{2\si (j)-1}}-x_{l_{2\si (j-1)-1}})\,
 u^{\ze'}(x_{l_{2\si' (j)-1}}-x_{l_{2\si' (j-1)-1}})\prod_{i=1}^{n}  \,dx_{l_{2i-1}}.\nn\\
  \label{2.0}
\end{eqnarray}
Writing $y_{\si (j)}=x_{l_{2\si (j)-1}}, y_{\si' (j)}=x_{l_{2\si' (j)-1}}$ and using (\ref{1.30g}) we can write this as
\begin{eqnarray}
&& 
\int  
\wt \mathcal{T}_{h}( x;\,\pi,\pi',e,e )\prod_{ j=1}^{ 2n}\,dx_{j}= \({8h^{3}(1+O(h))\over 3} \)^{n} \label{2.0a}\\
&&\hspace{
.3in}  \times\int \prod_{ j=1}^{ n}\,    u^{\ze}(y_{\si (j)}-y_{\si (j-1)})\,
u^{\ze'}(y_{\si' (j)}-y_{\si' (j-1)})\prod_{i=1}^{n}  \,dy_{i}.\nn
\end{eqnarray}

 Considering (\ref{4.28}),  (\ref{f9.41z2}) and (\ref{2.0a}) we see that  
 \begin{eqnarray} 
\lefteqn{ \int \mathcal{T}_{h}( x;\,\pi,\pi',e,e) \,\prod_{ j=1}^{ 2n}\,dx_{j}\label{1.0g}  =\({8h^{3}\over 3} \)^{n}\int \prod_{ j=1}^{ n}\,    u^{\ze}(y_{\si (j)}-y_{\si (j-1)})}\qquad\\
&&\hspace{1.2in}
\times u^{\ze'}(y_{\si' (j)}-y_{\si' (j-1)})\prod_{i=1}^{n}  \,dy_{i}+O(h^{3n+1}).\nn
 \end{eqnarray}
 
 In the first paragraph of  this subsection we explain  what we mean by $(\pi,\pi')\sim \mathcal{P}$,  for a pairing  $\mathcal{P}=\{(l_{2i-1},l_{2i})\,,\,1\leq i\leq n\}$  of the integers $[1,2n]$ and   permutations $\pi,\pi'$ of $[1,2n]$ that are  compatible with $\mathcal{P}$.  Obviously, there are many such pairs.   There are   $2^{2n}$ ways we can interchange the two elements of each pair $
\pi(2j-1), \pi(2j)$, and $
\pi'(2j-1), \pi'(2j)$  without changing (\ref{1.0g}).     Furthermore, by permuting the pairs $\{\pi(2j-1), \pi(2j)\}$ we give rise  to all possible permutations  $\si$ of $[1,n]$, and similarly for $\pi'$.   Consequently,
 \begin{eqnarray} 
\lefteqn{
  \sum_{(\pi,\pi')\sim \mathcal{P}} \int \mathcal{T}_{h}( x;\,\pi,\pi',e,e) \,\prod_{ j=1}^{ 2n}\,dx_{j} \nn}\\
&&  =\({32h^{3}\over 3} \)^{n}\sum_{\si,\si'}\int \prod_{ j=1}^{ n}\,    u^{\ze}(y_{\si (j)}-y_{\si (j-1)})\,
u^{\ze'}(y_{\si' (j)}-y_{\si' (j-1)})\prod_{i=1}^{n}  \,dy_{i}\nn\\
&&\hspace{3.5in}  +O(h^{3n+1}) \nn\\ 
&& =\({32h^{3}\over 3} \)^{n}  E\lc\(\int  L^{ x}_{ \la_{\ze}}  \wt L^{ x}_{\wt \la_{\ze'}}    \,dx\)^{
n}\rc+O(h^{3n+1}). \label{1.0h}
 \end{eqnarray} 
Here the  sum  in the second line runs over all   permutations  $\si,\,\si'$ of $\{1,\ldots, n\}$ and    $\si(0)=\si'(0)=0$. The      final   line of (\ref{1.0h}) follows from  the Kac moment formula,  (\ref{kac.m}).  

Since   there are $  ( 2n)!/( 2^{ n}n!)$ pairings of the $2n$
elements  $\{1,\ldots, 2n\}$ we obtain (\ref{f9.45sj}).  

\medskip  In the next two subsections we obtain (\ref{big.2}).

 \subsection{${\bf a=a'=e}$ without  compatible  permutations}\label{ss-3.2}
 
  Consider the multigraph $G_{\pi,\pi'}$ whose vertices consist of   
$\{1,\ldots, 2n\}$ and  assign an edge 
between the vertices $\pi (2j-1) $ and $ \pi (2j)$ for each $j=1,\ldots,n$ and similarly between $\pi' (2j-1) $ and $ \pi' (2j)$   for each $j=1,\ldots,n$. Each vertex is connected to two edges, and it is possible to have two edges between any two vertices $i,j$.   Note that the connected components $C_{j}$, $j=1,\ldots, k$ of $G_{\pi,\pi'}$ consist of cycles. (For example,  in   Subsection \ref{ss3.1}, all the cycles  are  of order two.)

 Let $C_{j}=\{j_{1},\ldots, j_{ l(j)}\}$ be written in cyclic order where $l(j)=|C_{j} |$.   Clearly
$\sum_{j=1}^{k}l(j)=2n$.
We show that when all the cycles are not of order two, as they are in the case of    compatible permutations considered in Subsection \ref{ss3.1},  then  
\begin{equation}
\bigg|\int \mathcal{T}_{h}( x;\,\pi,\pi',e,e) \,\prod_{ j=1}^{ 2n}\,dx_{j}\bigg|\leq C  h^{3n+1}.\label{f9.50}
\end{equation}

Since we    only need an upper bound, we  take absolute values in the integrand     and get  
\begin{eqnarray}
\lefteqn{\bigg|\int \mathcal{T}_{h}( x;\,\pi,\pi',e,e) \,\prod_{ j=1}^{ 2n}\,dx_{j}\bigg|
 \label{f9.51}}\\
&&   \leq   \int    \prod_{j=1}^{k}\( w(x_{j_{2}}-x_{j_{1}})\cdots    w(x_{j_{l(j)}}-x_{j_{l(j)-1}})w(x_{j_{1}}-x_{j_{l(j)}})   \) \nn\\ &&\hspace{
.5in}\times\prod_{ j=1}^{ n}\, u (x_{\pi(2j-1)}-x_{\pi(2j-2)})\,u(x_{\pi'(2j-1)}-x_{\pi'(2j-2)})   \, \prod_{ j=1}^{ 2n}\,dx_{j}\nn,
\end{eqnarray}
where we use the notation $u(x)$ to denote either $u^{\ze}(x)$ or $u^{\ze'}(x)$, and 
$w(x)$ to denote either $w^{\ze}(x)$ or $w^{\ze'}(x)$.  ($w^{\ze}(x)$ is defined in (\ref{4.20}).) 
  Note that we group the functions $w$ according to the cycles.

  For each $j=1,\ldots,k$  we set $y_{j_{i}}=x_{j_{i}}-x_{j_{i-1}}$, $i=2,\ldots,l(j),,$ and note that  $\sum_{i=2}^{l(j)}y_{j_{i}}=-(x_{j_{1}}-x_{j_{l(j)}})$.   It is easy to see that the $2n-k$ variables $\{y_{j_{i}}\,|\,j=1,\ldots,k\,; i=2,\ldots, l(j)\}$ are linearly independent. We  then choose an additional $k$ variables $z_{l}\,;\,l=1,\ldots,  k$
from  amongst the variables $\{x_{\pi(2j-1)}-x_{\pi(2j-2)}\,,\,x_{\pi'(2j-1)}-x_{\pi'(2j-2)}\,;1\leq j\leq n\}$ so that $\{y_{j_{i}}\,|\,j=1,\ldots,k\,; i=2,\ldots, l(j)\}\cup \{z_{l}\,|\,l=1,\ldots,  k\}$ are linearly  independent and generate $\{x_{1},\ldots, x_{2n}\}$.   We make this change of  variables and use the fact that $u(x)$ is bounded and integrable,   followed by   (\ref{1.3yq})  and (\ref{li.13}), to see that
\begin{eqnarray}
&&\bigg|\int \mathcal{T}_{h}( x;\,\pi,\pi',e,e) \,\prod_{ j=1}^{ 2n}\,dx_{j}\bigg|
\label{f9.53}\\
&&  \qquad \leq C   \prod_{j=1}^{k}\(\int   w(y_{j_{2}})\cdots    w(y_{j_{l(j)}} )w\(\sum_{i=2}^{l(j)}y_{j_{i}}\)    \,\,  \prod_{ i=2}^{ l(j)}\,dy_{j_{i}}\)\nonumber\\
&&  \qquad \leq C  \prod_{j=1}^{k}\, \sup_{x}|w(x)|\(\int   w(y )\,dy \)^{l(j)-1}\nonumber\\
&&\qquad \leq C \prod_{j=1}^{k}h^{1+2(l(j)-1)}=C \prod_{j=1}^{k}h^{2l(j)-1}.\nonumber
\end{eqnarray}
  (Note that the only dependence on $\ze$ and $\ze'$ is in the constant $C$.)

Since   $\sum_{j=1}^{k}l(j)=2n$, we see from (\ref{f9.53}) that  
\begin{eqnarray}
&&\bigg|\int \mathcal{T}_{h}( x;\,\pi,\pi',e,e) \,\prod_{ j=1}^{ 2n}\,dx_{j}\bigg|\leq C h^{4n-k}=
C h^{3n }h^{n-k}.
\label{f9.55}
\end{eqnarray}

It is easily seen that for non-compatible permutations we have $k<n$, which proves (\ref{f9.50}).

\subsection{When ${\bf a=a'=e}$  does  not hold}\label{ss3.4}
 
   We now consider all partitions   $\pi$ and $\pi'$ when $a=a'=e$ does not hold.  
  Consider the basic formula    (\ref{f1.21gi}).   Since we    only need an upper bound, we  take absolute values in the integrand as in (\ref{f9.51}).  Since ${  a=a'=e}$  does  not hold
 there are terms in which only one $\De^{h}$ is applied to a     $u^{\ze} $ or $u^{\ze'} $.  
 
We use the notation $u$ and $w$ defined right after (\ref{f9.51}).  If there are $k<2n$   factors of the  type  $w $, then there are  
 $2(2n-k)$ factors of the type    $\De^{\pm h}u$.   We use  (\ref{1.3x}) to   pull out a factor of 
 \begin{equation}
   h^{2(2n-k)}.\label{6.1}
 \end{equation}
 from  the basic formula   (\ref{f1.21gi}),
and are left with an integral like the one on the right--hand side of  (\ref{f9.51}), 
 except that there are $k$   factors of the form  $w $ which may be linked in chains as well as in cycles   and there   are $4n-k$ factors of type $u$. We denote this integral by $J_{h}$.     
 
   As  in (\ref{f9.51}), we arrange the $w$ factors into cycles and chains.
  We   then change variables and integrate   the $w$ factors. As  in  (\ref{f9.53}) a cycle of length $l$ gives a contribution that is bounded by     $  C h^{1+2(l-1)}=C h^{2l-1} $.  In addition, 
   by (\ref{li.13}),   chain of length $l'$ gives a contribution that is bounded by   $ C h^{ 2l'}$. 
   
   If there are $j$ cycles of lengths $l(i),\,i=1,\ldots,j$ and 
$j'$ chains of lengths $l'(i),\,i=1,\ldots,j'$, we  have
\begin{equation}
\sum_{i=1}^{j}l(i)+\sum_{i=1}^{j'}l'(i)=k.\label{f9.50cb}
\end{equation}
Therefore  
\bea
J_{h}&\leq &Ch^{\(2\sum_{i=1}^{j}l(i)\)-j} h^{2\sum_{i=1}^{j'}l'(i)}\nn\\
&\leq &C h^{2k-j} \label{f9.50cc}.
\eea
Together with (\ref{6.1}) this shows that   
\begin{equation}
 |\int \mathcal{T}_{h}( x;\,\pi,\pi',a,a') \,\prod_{ j=1}^{ 2n}\,dx_{j}|\leq C h^{4n-j}=C h^{3n } 
 h^{n-j} \label{f9.50a}.
\end{equation}
As   in (\ref{f9.55}) we see  that 
\begin{equation}
\Bigg |\int \mathcal{T}_{h}( x;\,\pi,\pi',a,a') \,\prod_{ j=1}^{ 2n}\,dx_{j}\Bigg|\leq C h^{3n+1}. \label{f9.50aa}
\end{equation}

\medskip  	We have  established (\ref{f9.45sj})   when $m$ is even. We  now show  that we get the same estimates when  $u^{\ze}$ and $u^{\ze'}$   are replaced by $u^{\ze,\sharp}$
  and   $u^{\ze',\sharp}$; (see   (\ref{3.9}) and (\ref{f1.21gi})).  
  
 The key observation that explains this is that in applying   the product formula (\ref{prf}), the only  terms of the form $u^{\ze}(x-y)$ that may have $x$ replaced by $x+ h$   are those   to which $\De^{h}_{x}$ is not applied.  Similarly $y$ may be replaced by $y+h$  only if $\De^{h}_{y}$ is not applied to a term of the form $u^{\ze}(x-y)$.   Consequently,  in evaluating (\ref{3.11}) with   $\TT_{h}$ replaced by $\TT'_{h}$  we still have $\De^{h}\De^{-h}u^{\ze,\sharp}=\De^{h}\De^{-h}u^{\ze}$ and similarly for  $\De^{h}\De^{-h}u^{\ze',\sharp}$. 

\label{page24}It is easy to see that the presence of the terms in $u^{\ze.\sharp}$  or $\De^{\pm h}u^{\ze,\sharp}$, or in $u^{\ze'.\sharp}$  or $\De^{\pm h}u^{\ze',\sharp}$    have no effect on the  integrals that are $O(h^{3n+1/2})$ as $h\to 0$.   (I.e. the terms that   are equal to 0 in (\ref{3.11}).)  This is because in evaluating these expressions we either integrate over all of $R^{1}$ or else use bounds that hold on all of  $R^{1}$.

 They do have an effect on the terms for which  the limit in (\ref{3.11}) are not zero. For example, instead of the right-hand side of (\ref{1.0g}), we now have

 \bea && \label{1.0gaa}   \({8h^{3}\over 3} \)^{n}\int \prod_{ j=1}^{ n}\,    u^{\ze,\sharp}(y_{\si (j)}-y_{\si (j-1)})  \\
 &&\quad\qquad
\times \,u^{\ze',\sharp}(y_{\si' (j)}-y_{\si' (j-1)})\prod_{i=1}^{n}  \,dy_{i}+O(h^{3n+1}).\nn
 \eea
  Suppose that $u^{\ze,\sharp}_{s_{i}}(y_{\si (i)}-y_{\si (i-1)})=u^{\ze } (y_{\si (i)}-y_{\si (i-1)}\pm h)$.  We write this term as
\begin{equation}\qquad
u^{\ze,\sharp}(y_{\si (i)}-y_{\si (i-1)})=u^{\ze }(y_{\si (i)}-y_{\si (i-1)}) +\De^{\pm h}u^{\ze }(y_{\si (i)}-y_{\si (i-1)}) \label{4.52}
   \end{equation}
   and similarly for $u^{\ze',\sharp}$.  Substituting these expressions  into  (\ref{1.0gaa})
 and using  (\ref{1.3x}) it is easy to see that (\ref{1.0gaa}) is asymptotically equivalent to  the right-hand side of (\ref{f9.45sj})   when $m$ is even. (The error term may be different). Thus we see that replacing  $u^{\ze}$ and $u^{\ze'}$      by $u^{\ze,\sharp}$
  and   $u^{\ze',\sharp}$    does not change (\ref{3.11}) when $m$ is even.
  
     \medskip	  It is rather obvious that the limit in (\ref{3.11}) is zero when $m$ is odd because in this case we can not construct a graph with all cycles of order 2. The extension of this limit when $u^{\ze}$ and $u^{\ze'}$   are replaced  by $u^{\ze,\sharp}$
  and   $u^{\ze',\sharp}$ follows as above.\qed

\section{Proof of Lemma \ref{lem-Lap}}\label{sec-Lap}
 
For $h=0$ it suffices to show that 
\begin{eqnarray} &&
G_{0}(s,t):= E\lc\(  \int  L^{ x}_{s}  \wt L^{ x}_{t}    \,dx \,\,\)^{
n}\rc \label{67.1e}
\end{eqnarray}
 is a non-negative, polynomially bounded continuous   increasing function of $(s,t)$. The fact that $G_{0}(s,t)$ is
 a non-negative,    increasing function of $(s,t)$ follows immediately from the fact that the local times $ L^{ x}_{s}$ and $ \wt L^{ x}_{t}  $ have these properties. 
 
 To prove continuity we note that for all $|r|,|r'|\le r_{0}$, $ \int  L^{ x}_{s+r}  \wt L^{ x}_{t+r'}    \,dx
 \le \int  L^{ x}_{s+r_{0}}  \wt L^{ x}_{t+r_{0}}    \,dx$. Therefore continuity follows from the Dominated Convergence Theorem and the continuity of local times once we show that for all $s,t$
 \begin{equation}
    \int  L^{ x}_{s}  \wt L^{ x}_{t}    \,dx\label{5.2}
   \end{equation}
 has all   moments.
 It follows from the Cauchy--Schwarz inequality,   the scaling  relationship (\ref{scale}), and (\ref{57rb.1}), that
 \begin{eqnarray}
\lefteqn{E\lc\(  \int  L^{ x}_{s}  \wt L^{ x}_{t}    \,dx \,\,\)^{
n}\rc
 \label{cs5}}\\
 &&\leq    E\lc\(  \int  \(L^{ x}_{s} \)^{2}    \,dx \,\,\)^{
n/2}\rc   E\lc\(  \int   \( \wt L^{ x}_{t} \)^{2}   \,dx \,\,\)^{
n/2}\rc \nonumber\\
 &&\leq  (st)^{3n/4}  E\lc\(  \int  \(L^{ x}_{1} \)^{2}    \,dx \,\,\)^{
n/2}\rc E\lc\(  \int   \( \wt L^{ x}_{1} \)^{2}   \,dx \,\,\)^{
n/2}\rc \nonumber\\
 &&\leq  C(st)^{3n/4}.  \nonumber
 \end{eqnarray}
In addition to showing that   (\ref{5.2}) has all moments, this also shows that $G_{0}(s,t)$   is a polynomially bounded function of $(s,t)$.

We now consider $F_{h}(s,t)$ for $h>0$. It suffices to show that 
\begin{eqnarray} &&
G_{h}(s,t):= E\(\(  \int ( L^{ x+h}_{s}- L^{ x}_{s})(\wt L^{ x+h}_{t}- \wt L^{ x}_{t})\,dx   \)^{m}\) \label{67.1e}
\end{eqnarray}
is a non-negative, polynomialy bounded, continuous   increasing function of $(s,t)$. 

Let $W_{t}$ denote Brownian motion and let $f\in \mathcal{S}(R^{1})$ be a positive symmetric function supported on $[-1,1]$ with $\int f(x)\,dx=1$. Set $f_{\ep}(x)=f(x/\ep)/\ep$ and
\begin{equation}
 L^{ x}_{s,\ep}=\int_{0}^{s}f_{\ep}(W_{r}-x)\,dr. \label{67.3}
\end{equation}
It follows from \cite[Lemma 2.4.1]{book} that
\begin{eqnarray}
&&E\(  \prod_{j=1}^{m} (L^{ x_{j}+h}_{s}- L^{ x_{j}}_{s})(\wt L^{ x_{j}+h}_{t}- \wt L^{ x_{j}}_{t}) \)
\label{par.1}\\
&&\qquad=\lim_{\ep\rar 0}E\(  \prod_{j=1}^{m} (L^{ x_{j}+h}_{s,\ep}- L^{ x_{j}}_{s,\ep})(\wt L^{ x_{j}+h}_{t,\ep}- \wt L^{ xj}_{t,\ep}) \)   \nonumber\\
&&\qquad=\lim_{\ep\rar 0}E\(  \prod_{j=1}^{m} (L^{ x_{j}+h}_{s,\ep}- L^{ x_{j}}_{s,\ep})  \) 
E\(  \prod_{j=1}^{m} (\wt L^{ x_{j}+h}_{t,\ep}- \wt L^{ x_{j}}_{t,\ep})  \).  \nonumber
\end{eqnarray}
Using the Fourier transform
\begin{equation}
L^{ x+h}_{s,\ep}- L^{ x}_{s,\ep}=\int e^{-ipx}(e^{-iph}-1)\wh f(\ep p)\int_{0}^{s}e^{ipW_{r}}\,dr\,dp\label{67.4}
\end{equation}
we have
\begin{eqnarray}
\lefteqn{ E\(  \prod_{j=1}^{m} (L^{ x_{j}+h}_{s,\ep}- L^{ x_{j}}_{s,\ep})  \)
\label{par.2}}\\
&& = \int_{R^{m}} \prod_{j=1}^{m}e^{-ip_{j}x_{j}}(e^{-ip_{j}h}-1)\wh f(\ep p_{j})\int_{[0,s]^{m}} E\(  \prod_{j=1}^{m}e^{ip_{j}W_{r_{j}}}  \)  \prod_{j=1}^{m}\,dr_{j}\,dp_{j}. \nonumber
\end{eqnarray}
Note that
\begin{eqnarray}
\lefteqn{\int_{[0,s]^{m}} E\(  \prod_{j=1}^{m}e^{ip_{j}W_{r_{j}}}  \) \prod_{j=1}^{m}\,dr_{j}
\label{par.4}}\\
&& =\sum_{\pi} \int_{\{0\leq r_{1}\leq \cdots\leq r_{m}\leq s \}} E\(  \prod_{j=1}^{m}e^{ip_{\pi(j)}W_{r_{j}}}  \) \prod_{j=1}^{m}\,dr_{j} \nonumber\\
&& =\sum_{\pi} \int_{\{0\leq r_{1}\leq \cdots\leq r_{m}\leq s \}} E\(  \prod_{j=1}^{m}e^{i(\sum_{k=j}^{m}p_{\pi(k)})(W_{r_{j}}-W_{r_{j-1}})}  \) \prod_{j=1}^{m}\,dr_{j} \nonumber\\
&& =\sum_{\pi} \int_{\{0\leq r_{1}\leq \cdots\leq r_{m}\leq s \}} \prod_{j=1}^{m}   e^{-(\sum_{k=j}^{m}p_{\pi(k)})^{2}( r_{j}- r_{j-1} )}   \prod_{j=1}^{m}\,dr_{j}. \nonumber
\end{eqnarray}
Since this is bounded and integrable in $p_{1},\ldots, p_{m}$, and $\wh f(\ep p)\le C$, we can take the limit as $\ep $ goes to zero in (\ref{par.2}) and hence in (\ref{par.1}), to see that
\begin{eqnarray}
\lefteqn{E\(  \prod_{j=1}^{m} (L^{ x_{j}+h}_{s}- L^{ x_{j}}_{s})(\wt L^{ x_{j}+h}_{t}- \wt L^{ x_{j}}_{t}) \)
\label{par.5}}\\
&& = \int_{R^{2m}} \prod_{j=1}^{m}e^{-i(p_{j}+p_{j}')x_{j}}(e^{-ip_{j}h}-1)(e^{-ip'_{j}h}-1) \nonumber\\
&&\qquad \times\int_{[0,s]^{m}} E\(  \prod_{j=1}^{m}e^{ip_{j}W_{r_{j}}}  \) \int_{[0,t]^{m}}E\(  \prod_{j=1}^{m}e^{ip'_{j}W_{r'_{j}}}  \) \prod_{j=1}^{m}\,dr'_{j}\,dr_{j}\,dp_{j}\,dp'_{j}. \nonumber
\end{eqnarray}
It now follows from Parseval's Theorem that
\begin{eqnarray}
\lefteqn{G_{h}(s,t)=\int  E\(  \prod_{j=1}^{m} (L^{ x_{j}+h}_{s}- L^{ x_{j}}_{s})(\wt L^{ x_{j}+h}_{t}- \wt L^{ x_{j}}_{t}) \)
 \prod_{j=1}^{m}\,dx_{j}\label{par.6}}\\
&& ={1 \over (2\pi)^{m}} \int_{R^{m}} \prod_{j=1}^{m} |e^{-ip_{j}h}-1|^{2} \nonumber\\
&&\qquad \times\int_{[0,s]^{m}} E\(  \prod_{j=1}^{m}e^{ip_{j}W_{r_{j}}}  \) \int_{[0,t]^{m}}E\(  \prod_{j=1}^{m}e^{ip_{j}W_{r'_{j}}}  \) \prod_{j=1}^{m}\,dr'_{j}\,dr_{j}\,dp_{j}. \nonumber
\end{eqnarray}
The fact that $G_{h}(s,t)$ is
 a non-negative,     increasing function of $(s,t)$   follows  from this and  (\ref{par.4}).
  The fact that $G_{h}(s,t)$ is
 a polynomialy bounded continuous function of $(s,t)$, follows as in the proof for $G_{0}(s,t)$ 
 if we note that by translation invariance $\int  \(L^{ x+h}_{s} \)^{2}    \,dx=\int  \(L^{ x}_{s} \)^{2}    \,dx$ so that, as in (\ref{cs5}),
  \begin{eqnarray}
\lefteqn{E\lc\(  \int  L^{ x+h}_{s}  \wt L^{ x}_{t}    \,dx \,\,\)^{
n}\rc
 \label{cs6}}\\
 &&\leq    E\lc\(  \int  \(L^{ x+h}_{s} \)^{2}    \,dx \,\,\)^{
n/2}\rc   E\lc\(  \int   \( \wt L^{ x}_{t} \)^{2}   \,dx \,\,\)^{
n/2}\rc \nonumber\\
 &&\leq    E\lc\(  \int  \(L^{ x}_{s} \)^{2}    \,dx \,\,\)^{
n/2}\rc   E\lc\(  \int   \( \wt L^{ x}_{t} \)^{2}   \,dx \,\,\)^{
n/2}\rc \nonumber\\
 &&\leq  C(st)^{3n/4}.  \nonumber
 \end{eqnarray}
 
 \qed

\section{Proof of Theorem \ref{theo-clt2}}\label{sec-BM}

The proof of  Theorem \ref{theo-clt2} follows from the next lemma exactly as in the proof of Theorem \ref{theo-cltindf} on page \pageref{page 5}.

\bl\label{lem-6.1} For each integer $m\ge 0$ and $t\in R_{+}$
\begin{eqnarray} &&
\lim_{ h\rar 0}E\(\({ \int ( L^{ x+h}_{t}- L^{ x}_{ t})^{
2}\,dx- 4h\over h^{ 3/2}}\)^{m}\)\nn\\ 
&&\hspace{ .5in}  =\left\{\begin{array}{ll}
\displaystyle{  ( 2n)!\over 2^{ n}n!}\( \displaystyle {64  \over 3}\)^{ n} E\lc\(\int (L^{ x}_{ t})^{ 2}\,dx\)^{
n}\rc &\mbox{ if }m=2n\\\\
0&\mbox{ otherwise.}
\end{array}
\right.
\label{7.53}
\end{eqnarray}
\el

We use the next lemma in the proof of Lemma \ref{lem-6.1}. It is proved in Section \ref{sec-clt}.  

\bl\label{lem-6.2}
Let $\la_{\ze}$ be an exponential random variable with mean $1/\ze$.
For each integer $m\ge 0$,
\begin{eqnarray} &&
\lim_{ h\rar 0}E\(\({ \int ( L^{ x+h}_{\la_{\ze}}- L^{ x}_{ \la_{\ze}})^{
2}\,dx-4h\la_{\ze}  \over h^{ 3/2}} \)^{m}\)\nn\\ 
&&\hspace{ .5in}  =\left\{\begin{array}{ll}
  \displaystyle  {( 2n)!\over 2^{ n}n!}\(   \displaystyle  {64 \over 3}\)^{ n} E\lc\(\int (L^{ x}_{ \la_{\ze}})^{ 2}\,dx\)^{
n}\rc &\mbox{ if }m=2n\\\\
0&\mbox{ otherwise.}
\end{array}
\right.
\label{7.54a}
\end{eqnarray}
  \el

\noindent{\bf Proof of Lemma \ref{lem-6.1} }We  write (\ref{7.54a}) as 
\begin{eqnarray} &&
\lim_{ h\rar 0}\int_{0}^{\ff}e^{- \ze s } E\(\({  \int ( L^{ x+h}_{s}- L^{ x}_{ s})^{
2}\,dx-4hs  \over h^{ 3/2}} \)^{m}\) \,ds 
\label{77.1d}\\
&&\qquad= \int_{0}^{\ff} e^{- \ze s }E\lc  \eta^{m}\( {64 \over 3}\int  (L^{ x}_{s})^{2}   \,dx\,\,\)^{
m/2}\rc \,ds. \nn
\end{eqnarray}
 Let  
\begin{eqnarray} &&
\wh F_{m,h}(s):= E\(\({\int ( L^{ x+h}_{s}- L^{ x}_{ s})^{2}\,dx-4hs 
 \over h^{ 3/2}} \)^{m}\),\hspace{.2 in}h>0 \label{77.1e}
\end{eqnarray}
and
\begin{eqnarray} &&
\wh F_{m,0}(s):= E\lc  \eta^{m}\( {64 \over 3}\int  (L^{ x}_{s})^{2}   \,dx\,\,\)^{
m/2}\rc. \label{77.1e}
\end{eqnarray}
Then (\ref{77.1d}) can be written as 
\begin{equation}
\lim_{ h\rar 0}\int_{0}^{\ff}e^{- \ze s } \wh {F}_{m,h}(s) \,ds =\int_{0}^{\ff}e^{- \ze s } \wh {F}_{m,0}(s) \,ds.\label{77.4}
\end{equation}

We consider first the case when $m$ is even and write   $m=2n$. In this case  $\wh {F}_{2n,h}(s)\geq 0$ and the extended continuity theorem 
\cite[XIII.1, Theorem 2a]{Feller} applied to  (\ref{77.4}) implies that
\begin{equation}
\lim_{ h\rar 0}\int_{0}^{t}  \wh {F}_{2n,h}(s) \,ds =\int_{0}^{t}  \wh {F}_{2n,0}(s) \,ds\label{77.5}
\end{equation}
for all $t$. In particular,
\begin{equation}
\lim_{ h\rar 0}\int_{t}^{t+\de}  \wh {F}_{2n,h}(s) \,ds =\int_{t}^{t+\de}  \wh {F}_{2n,0}(s) \,ds.\label{77.6}
\end{equation}
It follows from the proof of Lemma \ref{sec-Lap} that $ \wh {F}_{2n,0}(t)$ is continuous in $t$. Consequently,
\begin{equation}
\lim_{ \de\rar 0}\lim_{ h\rar 0}{1 \over \de}\int_{t}^{t+\de}  \wh {F}_{2n,h}(s) \,ds = \wh {F}_{2n,0}(t). \label{77.6}
\end{equation}
When $t=0$ we get
\begin{equation}
\lim_{\de\to 0^+}
\lim_{h\to 0^+}{1\over \de}\int_0^{\de}\wh {F}_{2n,h}(s)ds=0.\label{77.13}
\end{equation}

To obtain (\ref{7.53}) when $m$ is even we must show that
\begin{equation}
\lim_{ h\rar 0}   \wh {F}_{2n,h}(t)  = \wh {F}_{2n,0}(t) \label{77.15}.
\end{equation}
This follows from  (\ref{77.6}) once we  show that
\begin{equation}
\lim_{ \de\rar 0}\lim_{ h\rar 0}{1 \over \de}\int_{t}^{t+\de}  \wh {F}_{2n,h}(s) \,ds = \lim_{ h\rar 0}   \wh {F}_{2n,h}(t)  . \label{77.6j}
\end{equation}
We proceed to obtain (\ref{77.6j}).

For $s\geq t$ we  write 
\begin{eqnarray} \lefteqn{
\int ( L^{ x+h}_{s}- L^{ x}_{ s})^{2}dx-4hs
=\bigg\{\int ( L^{ x+h}_{t}- L^{ x}_{ t})^{2}dx-4ht\bigg\}\label{77.7}}\\
&&\qquad
+\bigg\{\int \Big[\big(L_{s}^{x+h}-L_{s}^x\big) -
\big(L_t^{x+h}-L_t^x\big)\Big]^2dx -4h(s-t)\bigg\}\nn\\
&&\qquad\quad
+2\bigg\{\int (L_t^{x+h}-L_t^x
)\Big[ (L_{s}^{x+h}-L_{s}^x ) -
\big(L_t^{x+h}-L_t^x\big)\Big]dx\bigg\}.\nn
\end{eqnarray}
We use the triangle inequality with respect to the norm $\|\,\cdot\,\|_{2n}$ to see that
\begin{eqnarray} \lefteqn{
\wh {F}_{2n, h}^{1/(2n)}(s)\le \wh {F}_{2n, h}^{1/(2n)}(t)
 \label{77.7aa}}\\
&&
+\bigg\{\E\bigg[{1\over h^{3/2}}\bigg\{\int \Big[\big(L_{s}^{x+h}-L_{s}^x\big) -
\big(L_t^{x+h}-L_t^x\big)\Big]^2dx -4h(s-t)\bigg\}\bigg]^{2n}\bigg\}^{1\over 2n}\nn\\
&&
+2\bigg\{\E\bigg[{1\over h^{3/2}}\bigg\{\int (L_t^{x+h}-L_t^x
)\Big[ (L_{s}^{x+h}-L_{s}^x ) -
\big(L_t^{x+h}-L_t^x\big)\Big]dx\bigg\}\bigg]^{2n}\bigg\}^{1\over 2n}.\nn
\end{eqnarray}
Note that  
\be \int \Big[\big(L_{s}^{x+h}-L_{s}^x\big) -
\big(L_t^{x+h}-L_t^x\big)\Big]^2dx\buildrel \LL\over =\int (L_{s-t}^{x+h}-L_{s-t}^x
)^2dx\label{77.9}
\ee
and
\begin{eqnarray} &&\int \big[L_t^{x+h}-L_t^x
\big]\Big[\big(L_{s}^{x+h}-L_{s}^x\big) -
\big(L_t^{x+h}-L_t^x\big)\Big]dx\label{77.8}\\
&&\qquad\buildrel \LL\over =\int 
\big[L_t^{x+h}-L_t^x\big]\big[\widetilde{L}_{s-t}^{x+h}-
\widetilde{L}_{s-t}^x\big]dx.\nn
\end{eqnarray}
Hence we can write (\ref{77.7aa}) as
\begin{eqnarray} 
\wh {F}_{2n, h}^{1/(2n)}(s)&\le &\wh {F}_{2n, h}^{1/(2n)}(t)
+\wh {F}_{2n, h}^{1/(2n)}(s-t)\label{77.7a}\\
&& 
+2\bigg\{\E\bigg[{1\over h^{3/2}}\int 
\big[L_{t}^{x+h}-L_{t}^x\big]\big[\widetilde{L}_{s-t}^{x+h}-
\widetilde{L}_{s-t}^x\big]dx\bigg]^{2n}\bigg\}^{1\over 2n}\nn.
\end{eqnarray}
We now use the triangle inequality   with respect to the norm in $L^{2n}([t,t+\de],   \de^{-1}\,ds)$ to see that
\begin{eqnarray} \lefteqn{
\bigg\{{1\over\delta}\int_t^{t+\delta}\wh {F}_{2n, h}(s)ds
\bigg\}^{1\over 2n}
\le \wh {F}_{2n, h}^{1/(2n)}(t)
+\bigg\{{1\over\delta}\int_t^{t+\delta}\wh {F}_{2n, h}(s-t)ds\bigg\}^{1\over 2n}\label{77.7b}}\qquad\\
&&+2\bigg\{{1\over\delta}\int_t^{t+\delta}\E\bigg[{1\over h^{3/2}}
\int 
\big[L_{t}^{x+h}-L_{t}^x\big]\big[\widetilde{L}_{s-t}^{x+h}-
\widetilde{L}_{s-t}^x\big]dx\bigg]^{2n}ds\bigg\}^{1\over 2n}\nn
\end{eqnarray}

A similar argument starting with (\ref{77.7}) shows that \begin{eqnarray} \lefteqn{
\bigg\{{1\over\delta}\int_t^{t+\delta}\wh {F}_{2n, h}(s)ds
\bigg\}^{1\over 2n}
\ge \wh {F}_{2n, h}^{1/(2n)}(t)
-\bigg\{{1\over\delta}\int_t^{t+\delta}\wh {F}_{2n, h}(s-t)ds\bigg\}^{1\over 2n}\label{77.7c}}\qquad\\
&&-2\bigg\{{1\over\delta}\int_t^{t+\delta}\E\bigg[{1\over h^{3/2}}
\int 
\big[L_{t}^{x+h}-L_{t}^x\big]\big[\widetilde{L}_{s-t}^{x+h}-
\widetilde{L}_{s-t}^x\big]dx\bigg]^{2n}ds\bigg\}^{1\over 2n}\nn
\end{eqnarray}
Since
\begin{equation}
{1\over\delta}\int_t^{t+\delta}\wh {F}_{2n, h}(s-t)ds={1\over\delta}\int_0^{\delta}\wh {F}_{2n, h}(s)ds\label{77.13n}
\end{equation}
we see from  (\ref{77.13}) that to prove (\ref{77.6j}) it only remains to show that
\begin{equation}\qquad
\lim_{\delta\to 0^+}\limsup_{h\to 0}
{1\over\delta}\int_t^{t+\delta}\E\bigg[{1\over h^{3/2}}
\int 
\big[L_{t}^{x+h}-L_{t}^x\big]\big[\widetilde{L}_{s-t}^{x+h}-
\widetilde{L}_{s-t}^x\big]dx\bigg]^{2n}ds=0.
\label{77.7d}
\end{equation}

By the monotonicity property of $F_{h}(s,t;m)$  given in Lemma \ref{lem-Lap},
\begin{eqnarray} &&
{1\over\delta}\int_t^{t+\delta}\E\bigg[{1\over h^{3/2}}
\int 
\big[L_{t}^{x+h}-L_{t}^x\big]\big[\widetilde{L}_{s-t}^{x+h}-
\widetilde{L}_{s-t}^x\big]dx\bigg]^{2n}ds\label{77.7f}\\
&&\qquad\le \E\bigg[{1\over h^{3/2}}
\int_{-\infty}^\infty
\big[L_{t}^{x+h}-L_{t}^x\big]\big[\widetilde{L}_\delta^{x+h}-
\widetilde{L}_\delta^x\big]dx\bigg]^{2n}\nn
\end{eqnarray}
Thus (\ref{77.7d}) follows from the fact that
\begin{equation}
\lim_{\delta\to 0^+}\limsup_{h\to 0^+}\E\bigg[{1\over h^{3/2}}
\int 
(L_t^{x+h}-L_t^x)(\widetilde{L}_\delta^{x+h}-
\widetilde{L}_\delta^x)dx\bigg]^{2n}=0\label{77.12}
\end{equation}
which, itself,  is a simple consequence of (\ref{57.4}) and Lemma \ref{lem-Lap}. Thus we obtain (\ref{77.15}) and hence  (\ref{7.53}) when $m$ is even.

\medskip	
In order to obtain (\ref{77.15}) when  $m$ is odd we first show that 
\begin{equation}
\sup_{h>0}\wh F_{2n,h}(t)\leq Ct^{3n/2}.\label{eq.1}
\end{equation}
To see this we  observe that by first changing variables and then using  the scaling relationship (\ref{scale}) with $h=\sqrt{t}$,
we have 
\begin{eqnarray}
\int ( L^{ x+h}_{t}- L^{ x}_{ t})^{2}\,dx&=&\sqrt{t}\int ( L^{\sqrt{t}( x+ht^{-1/2})}_{t}- L^{\sqrt{t} x}_{ t})^{2}\,dx\label{eq.1a}\\
&=& t^{3/2}\int ( L^{  x+ht^{-1/2}}_{1}- L^{  x}_{ 1})^{2}\,dx\nonumber.
\end{eqnarray}
Therefore
\begin{eqnarray}\qquad
{\int ( L^{ x+h}_{t}- L^{ x}_{ t})^{2}\,dx-4ht 
 \over h^{ 3/2}}&\stackrel{\mathcal{L}}{=}& \label{eq.2} 
  {t^{3/2}\(\int ( L^{  x+ht^{-1/2}}_{1}- L^{  x}_{ 1})^{2}\,dx-4ht^{-1/2}\) 
 \over h^{ 3/2}}\nn \\
&=& t^{3/4}\,\,{\(\int ( L^{  x+ht^{-1/2}}_{1}- L^{  x}_{ 1})^{2}\,dx-4ht^{-1/2}\) 
 \over (ht^{-1/2})^{ 3/2}}  \nonumber\\
\end{eqnarray}
so that for any integer $m$
\begin{equation}
\wh F_{m,h}(t)=t^{3m/4}\wh F_{m,ht^{-1/2}}(1).\label{eq.3}
\end{equation}
Therefore to prove (\ref{eq.1}) we  need only show that
\begin{equation}
\sup_{t}\sup_{ h>0 }\wh F_{2n,ht^{-1/2}}(1)\leq C.\label{eq.4}
\end{equation}

It follows from (\ref{77.15}) that for some $\de>0$
\begin{equation}
\sup_{\{t,h\,|\,ht^{-1/2}\leq \de\}} \wh{F}_{2n,ht^{-1/2}}(1)\leq C.\label{eq.5}
\end{equation}
 On the other hand, for $ht^{-1/2}\geq \de$
\begin{eqnarray}
&&\Bigg|{\(\int ( L^{  x+ht^{-1/2}}_{1}- L^{  x}_{ 1})^{2}\,dx-4ht^{-1/2}\) 
 \over (ht^{-1/2})^{ 3/2}}\Bigg|
\label{eq.6}\\
&& \qquad\leq \de^{-3/2}\int ( L^{  x+ht^{-1/2}}_{1}- L^{  x}_{ 1})^{2}\,dx +4\de^{-1/2} \nonumber\\
&&\qquad \leq 4\de^{-3/2}\int (L^{  x}_{ 1})^{2}\,dx +4\de^{-1/2}<\ff \nonumber
\end{eqnarray}
  since  $\int (L^{  x}_{ 1})^{2}\,dx$ has finite moments.  (See  (\ref{57rb.1}). Using (\ref{eq.5}) and (\ref{eq.6}) we get (\ref{eq.4}) and hence (\ref{eq.1}). It
then follows from   the Cauchy-Schwarz inequality that
\begin{equation}
\sup_{h>0}|\wh{F}_{m,h}(t)|\leq Ct^{3m/4}\label{eq.7}
\end{equation}
 for all integers $m$.
 
We next show that for any integer  $m$, the family of functions $\{\wh{F}_{m,h}(t);\,h\}$ is  equicontinuous in $t$, that is, for each $t$ and $\ep>0$ we can find a $\de>0$ such that  
\begin{equation}
\sup_{\{s\,|\,|s-t|\leq \de \}}\sup_{ h>0 }|\wh{F}_{m,h}(t)-\wh{F}_{m,h}(s)  |\leq \ep.\label{eq.8}
\end{equation}
Let
\begin{equation}
\Phi_{h}(t):={\int ( L^{ x+h}_{t}- L^{ x}_{ t})^{2}\,dx-4ht 
 \over h^{ 3/2}}\label{eq.9}
\end{equation}
Applying the identity
$A^{m}-B^{m}=\sum_{j=0}^{m-1}A^{j}(A-B)B^{m-j-1}$  with $A=\Phi_{h}(t),\,B=\Phi_{h}(s)$ gives
\begin{equation}
   \wh{F}_{m,h}(t)-\wh{F}_{m,h}(s)=\sum_{j=0}^{m-1}\Phi_{h}(t)^{j}(\Phi_{h}(t)-\Phi_{h}(s))\Phi_{h}(s)^{m-j-1}
   \end{equation}
Consequently by using  the Cauchy--Schwarz inequality twice and  (\ref{eq.7}),  we see that
\begin{equation}
\sup_{\{s\,|\,|s-t|\leq \de \}}\sup_{ h>0 }|\wh{F}_{m,h}(t)-\wh{F}_{m,h}(s)  |\leq C_{t,m}\sup_{\{s\,|\,|s-t|\leq \de \}}\sup_{ h>0 }\| \Phi_{h}(t)-\Phi_{h}(s)  \|_{2}.\label{eq.10}
\end{equation}
Using (\ref{77.7})--(\ref{77.8}), we see that to obtain (\ref{eq.8}) it  suffices to show that for some $\de>0$
\begin{equation}
\sup_{\{s\,|\,s\leq \de \}}\sup_{ h>0 }\wh{F}_{2,h}(s)  \leq \ep\label{eq.11}
\end{equation}
and  for any $T<\ff$
\begin{equation}
\sup_{\{\,t\leq T \}}\sup_{\{ \,s\leq \de \}}\sup_{ h>0 }\E\bigg[{1\over h^{3/2}}
\int 
(L_t^{x+h}-L_t^x)(\widetilde{L}_s^{x+h}-
\widetilde{L}_s^x)\,dx\bigg]^{2}\leq \ep.\label{eq.12}
\end{equation}

By (\ref{eq.1})
\begin{equation}
\sup_{h>0}F_{2,h}(s)\leq Cs^{3/2}\label{eq.13},
\end{equation}
which immediately gives (\ref{eq.11}).  Furthermore,  applying the Cauchy--Schwarz inequality in (\ref{par.6}) and using  (\ref{par.4}) to see that 
\begin{equation}
   \int_{[0,t]^{m}} E\(  \prod_{j=1}^{m}e^{ip_{j}W_{r_{j}}}  \) \prod_{j=1}^{m}\,dr_{j}
   \end{equation}
is positive and increasing in $t$,
 we see     that for  all $t\leq T$
\begin{eqnarray}
&&\E\bigg[{1\over h^{3/2}}
\int 
(L_t^{x+h}-L_t^x)(\widetilde{L}_s^{x+h}-
\widetilde{L}_s^x)\,dx\bigg]^{2}
\label{eq.14}\\
&&\qquad\leq \(\E\bigg[{1\over h^{3/2}}
\int 
(L_t^{x+h}-L_t^x)(\widetilde{L}_t^{x+h}-
\widetilde{L}_t^x)\,dx\bigg]^{2}\)^{1/2}   \nonumber\\
&&\qquad\qquad\times  \(\E\bigg[{1\over h^{3/2}}
\int 
(L_s^{x+h}-L_s^x)(\widetilde{L}_s^{x+h}-
\widetilde{L}_s^x)\,dx\bigg]^{2}\)^{1/2}   \nonumber\\
&&\qquad\leq \(\E\bigg[{1\over h^{3/2}}
\int 
(L_T^{x+h}-L_T^x)(\widetilde{L}_T^{x+h}-
\widetilde{L}_T^x)\,dx\bigg]^{2}\)^{1/2}   \nonumber\\
&&\qquad\qquad\times  \(\E\bigg[{1\over h^{3/2}}
\int 
(L_s^{x+h}-L_s^x)(\widetilde{L}_s^{x+h}-
\widetilde{L}_s^x)\,dx\bigg]^{2}\)^{1/2}   \nonumber.
\end{eqnarray}
Using the scaling relationship,   as in   (\ref{eq.2}), we see that
\begin{eqnarray}
&&
\E\bigg[{1\over h^{3/2}}
\int 
(L_s^{x+h}-L_s^x)(\widetilde{L}_s^{x+h}-
\widetilde{L}_s^x)\,dx\bigg]^{2}\label{eq.15}\\
&&\qquad  =s^{3/4} \E\bigg[{1\over (hs^{-1/2})^{3/2}}
\int 
(L_1^{x+hs^{-1/2}}-L_1^x)(\widetilde{L}_1^{x+hs^{-1/2}}-
\widetilde{L}_1^x)\,dx\bigg]^{2}.\nn
\end{eqnarray}
Following the proof of (\ref{eq.4}) we see that the expectation is bounded in $s$ and $h$. Therefore, by taking $\de$ sufficiently small we get (\ref{eq.12}).
This establishes (\ref{eq.8}).

\medskip	 We now obtain (\ref{7.53}) when $m$ is odd. By equicontinuity, for any sequence $h_{n}\rar 0$, we can find a 
 subsequence $h_{n_{j}}\rar 0$,  such that 
 \begin{equation}
\lim_{j\rar\ff}\wh F_{m,h_{n_{j}}}(t)\label{eq.16}
 \end{equation}
converges   to a continuous function which we denote by $\ov F_{m}(t)$. It remains to show that
\begin{equation}
\ov F_{m}(t)\equiv 0.\label{eq.17}
\end{equation}
Let 
\begin{equation}
G_{m,h }(t):=e^{-t}\wh F_{m,h }(t) \hspace{.2 in}\mbox{and}\qquad\ov G_{m}(t):=e^{-t}\ov F_{m}(t).\label{eq.18}
\end{equation}
By (\ref{eq.7})
\begin{equation}
\sup_{h>0}\sup_{t}|G_{m,h}(t)|\leq C\quad\mbox{ and}\quad\lim_{t\rar \ff}\sup_{h>0}G_{m,h}(t)=0. \label{eq.19}
\end{equation}
 It then  follows from (\ref{77.4}) and the dominated convergence theorem that for all $\ze>0$
\begin{equation}
\int_{0}^{\ff}e^{- \ze s } \ov G_{m}(s) \,ds =0.\label{eq.20}
\end{equation}
We obtain  (\ref{eq.17}) by showing that $\ov G_{m}(s)\equiv 0$.

It   follows from (\ref{eq.19}) that  $\ov G_{m}(t) $ is a continuous bounded function on $R_{+}$ that goes to zero as $t\to\ff$. By the Stone--Weierstrass Theorem; (see \cite[Lemma 5.4]{K}), for any $\ep>0$, we can find a finite linear combination of the form $\sum_{i=1}^{n}c_{i}e^{- \ze_{i} s }$ such that
\begin{equation}
\sup_{t}|\ov G_{m}(t)-\sum_{i=1}^{n}c_{i}e^{- \ze_{i} t }|\leq \ep.\label{eq.21}
\end{equation}
Therefore, by  (\ref{eq.20})
\begin{eqnarray} 
\int_{0}^{\ff}e^{-  s } \ov G^{2}_{m}(s) \,ds\label{eq.22}&
 =&\int_{0}^{\ff}e^{-  s }\(\sum_{i=1}^{n}c_{i}e^{- \ze_{i} s }\) \ov G_{m}(s) \,ds\\
 &&\quad+\int_{0}^{\ff}e^{-  s } \(\ov G_{m}(s)-\sum_{i=1}^{n}c_{i}e^{- \ze_{i} s }\)\ov G_{m}(s) \,ds\nn\\
 &=&\int_{0}^{\ff}e^{-  s } \(\ov G_{m}(s)-\sum_{i=1}^{n}c_{i}e^{- \ze_{i} s }\)\ov G_{m}(s) \,ds\nn\\
 &\le &2\ep\(\int_{0}^{\ff}e^{-  s } \ov G^{2}_{m}(s) \,ds\)^{1/2}
\end{eqnarray}
by the Cauchy--Schwarz inequality and  (\ref{eq.21}).   Thus $\int_{0}^{\ff}e^{-  s } \ov G^{2}_{m}(s) \,ds=0$ which implies that $\ov G_{m}(s)\equiv 0$.
\qed

\section{\bf Proof of Lemma   \ref{lem-6.2}}\label{sec-clt} 

\noindent{\bf  Proof of Lemma \ref{lem-6.2}$\,$}  
Our goal is the obtain  the asymptotic behavior  of the $m$--th moment of 
\begin{equation}
{ \int ( L^{ x+h}_{\la_{\ze}}- L^{ x}_{ \la_{\ze}})^{
2}\,dx-4h\la_{\ze}  \over h^{ 3/2}}\label{mot.1}
\end{equation}
as $h\rar 0$.  In the numerator we     have the term  $4h\la_{\ze}$. Note that  by Lemma \ref{lem-vare}, this is necessary in order  that the expected value of the numerator goes to $0$.  Since we have $h^{3/2}$ in the denominator in (\ref{mot.1}),  and  $O(h/h^{3/2})=O(h^{-1/2})$, we must show that in the expansion of the expectation of the $m$--th   moment of (\ref{mot.1}), the terms that would cause it to blow up  are canceled. We do this   in the first part of this proof.

 \medskip	
Note that  
\begin{equation}
   \int  L^{ x}_{\la_{\ze}}\,dx= \la_{\ze}.\label{7.1}
   \end{equation}
Using this and (\ref{1.8}), we write the left--hand side of (\ref{7.54a}) as 
\be  
\lim_{ h\rar 0}E\(\({ \int ( L^{ x+h}_{\la_{\ze}}- L^{ x}_{ \la_{\ze}})^{
2}\,dx- 2\De^{h}\De^{-h}u^{\ze}(0)\int  L^{ x}_{\la_{\ze}}\,dx\over h^{ 3/2}}\)^{m}\).\label{7.2}
 \ee
For any
integer
$m$ we have  
\bea \lefteqn{ E\(\( \int ( L^{ x+h}_{ \la_{\ze}}- L^{ x}_{ \la_{\ze}})^{ 2}\,dx-2\De^{h}\De^{-h}u^{\ze}(0) \int  L^{ x}_{
\la_{\ze}}\,dx\)^{ m}
\)\label{1.16g}}\\ &&=E\(\prod_{ i=1}^{ m}\( \int ( L^{ x_{ i}+h}_{ \la_{\ze}}- L^{ x_{
i}}_{
\la_{\ze}})^{ 2}\,dx_{ i}-2\De^{h}\De^{-h}u^{\ze}(0)\int  L^{ x_{ i}}_{
\la_{\ze}}\,dx_{ i}\)
\)\nn\\ &&=\sum_{A\subseteq \{ 1,\ldots,m\}} ( -1)^{|A^{c}|}E\(
\(\prod_{ i\in A} \int ( L^{ x_{ i}+h}_{ \la_{\ze}}- L^{ x_{ i}}_{ \la_{\ze}})^{ 2}\,dx_{
i}\)\right.\nn\\ &&
\hspace{ 1.7in}\times\left.\(\prod_{ i\in A^{ c}}2\De^{h}\De^{-h}u^{\ze}(0)\int L^{ x_{ i}}_{
\la_{\ze}}\,dx_{ i}\) \).\nn
\eea

 We now show that there are many cancelations in the final equation in  (\ref{1.16g}), that eliminate  the problematical terms we discussed in the beginning of this proof,  and also  significantly simplifies  it.

 \medskip	Consider a generic term in the final equation in (\ref{1.16g}) without the integrals. 
To  clarify what is going on we 
   calculate 
\begin{equation} E\(\prod_{ i\in A} (L^{ x_{ i}+h}_{ \la_{\ze}}- L^{ x_{ i}}_{ \la_{\ze}}) (L^{
y_{ i}+h}_{ \la_{\ze}}- L^{ y_{ i}}_{ \la_{\ze}})
\prod_{ i\in A^{ c}}2\De^{h}\De^{-h}u^{\ze}(0)L^{ x_{ i}}_{\la_{\ze}}\),\label{1.18g}
\end{equation}
  keeping in mind that   $y_{ i}=x_{ i}$ for all $1\le i\le m$. Using the Kac moment formula, (\ref{1.2w}), we have  
  \begin{eqnarray} \lefteqn{ E\(\prod_{ i\in A}(L^{ x_{ i}+h}_{ \la_{\ze}}- L^{ x_{ i}}_{ \la_{\ze}})
(L^{ y_{ i}+h}_{ \la_{\ze}}- L^{ y_{ i}}_{ \la_{\ze}})\prod_{ i\in A^{ c}}2\De^{h}\De^{-h}u^{\ze}(0)L^{ x_{
i}}_{\la_{\ze}}\)\label{1.19g}}\\ && =\(\prod_{ i\in A}\De_{ x_{ i}}^{ h}\De_{ y_{ i}}^{
h}\)E\(\prod_{ i\in A}L^{ x_{ i}}_{
\la_{\ze}} L^{ y_{ i}}_{ \la_{\ze}}\prod_{ i\in A^{ c}}2\De^{h}\De^{-h}u^{\ze}(0)L^{ x_{ i}}_{\la_{\ze}}\)
\nonumber\\ && =\(2\De^{h}\De^{-h}u^{\ze}(0)\)^{ | A^{ c} |}\(\prod_{ i\in A}\De_{ x_{ i}}^{
h}\De_{ y_{ i}}^{ h}\)\,\,
\sum_{ \si\in \mathcal{B}_{A}}\prod_{ j=1}^{m+|A|}u^{\ze}( \si ( j)-\si ( j-1)) \nonumber
\end{eqnarray} where the sum runs over $\mathcal{B}_{A}$, the set of  all
bijections  
\be
\si :\,[1,\ldots,m+|A|]\mapsto 
\{ x_{ i}, y_{ i},  i\in A\}\cup \{ x_{ i}, i\in A^{ c} \}.\label{7.6}
\ee 

  As we did in the beginning of Section \ref{sec-expind} we use the product
rule 
\begin{equation}\De_{ x}^{ h}\{ f( x)g( x)\}=
\{ \De_{ x}^{ h} f( x)\}g( x+h)+f( x)\{  \De_{ x}^{ h}g( x)\}\label{pr}
\end{equation}
to expand the  last line  of (\ref{1.19g}) into a sum of many terms, over all   
$\si\in \mathcal{B}_{A}$ and all ways to allocate each difference operator, $\De_{ x_{ i}}^{ h}$ and    $\De_{ y_{ j}}^{ h}$, $  i,j\in A$, to
the terms $u^{\ze}( \si ( j)-\si ( j-1))$ in which $\si (j-1)$ and/or $\si (j)$ are contained in $A$. After setting all $y_{ i}=x_{ i}$ we can then write (\ref{1.19g}) as
\begin{eqnarray}
&&\(2\De^{h}\De^{-h}u^{\ze}(0)\)^{ | A^{ c} |}
\label{nm.1}\\
&&  \sum_{ \si\in \mathcal{B}_{A},\,a} \prod_{ j=1}^{ m+|A|}\(\De^{ h}_{ \si(j)}\)^{a_{ 1}(j)}
\(\De^{ h}_{ \si(j-1)}\)^{a_{ 2}(j)}\,u^{\ze,\sharp}(  \si(j)- \si(j-1))\,\,\Big |_{y_{ i}=x_{ i},\,\forall i  } \nonumber
\end{eqnarray}
where   the sum runs over   $ \si\in \mathcal{B}_{A}$  and   all   
$a =(a_{ 1},a_{ 2})\,:\,[1,\ldots, m+|A|]\mapsto \{ 0,1\} \times \{ 0,1\} $, with the
restriction that for each $i\in A$ there is exactly one  factor  of the form $\De^{
h}_{ x_{i}}$ and one  factor  of the form $\De^{
h}_{ y_{i}}$, and there are no such factors for $i\in A^{c}$.  (Here we define $(\De_{x_{i}}^{h})^{0}=1 $ and $(\De_{0}^{h}) =1 $.) In this formula, $u^{\ze,\sharp}(x)$ can take any of the values $u^{\ze}(x)$, $u^{\ze}(x+h)$ or $u^{\ze}(x-h)$.  (This is because we use  (\ref{pr}) to pass from the last line of  (\ref{1.19g}) to (\ref{nm.1}).  We consider all three possibilities in the subsequent proofs.) It is important to recognize that in  (\ref{nm.1}) each of  the difference operators is applied to only one of the  terms   $u^{\ze,\sharp}(\cdot)$.

We get the simplification of the final equation in (\ref{1.16g}), because    many terms in the expansion of (\ref{1.19g}) for different sets $A$ and  $\si\in \mathcal{B}_{A}$  are the same, and when they are added, as they are in the final equation in (\ref{1.16g}), they cancel. We now make this precise.

Fix $A\subseteq \{ 1,\ldots,m\}$ and consider a particular bijection $\si\in \mathcal{B}_{A}$. Consider  (\ref{nm.1})
for this $A$ and $\si $. For $i\in A$ we
  say that
$x_{ i}$ is a  bound variable,   if $x_{ i}$ and $y_{ i} $ are adjacent,    i.e., if either $(x_{ i},y_{ i})=(\si (j-1),\si (j))$ or $(y_{ i},x_{ i})=(\si (j-1),\si (j))$ for some $j$.    Furthermore, for a given   $\si\in\BB_{A} $ that contains bound variables,  and a given $a$,
we say that a bound variable $x_{ i}$ is a  singular variable  
if both $\De_{ x_{ i}}^{ h}$ and
$\De_{ y_{ i}}^{ h}$ are  applied to the factor $u^{\ze}( x_{ i}-y_{ i})$.  

  Note that by (\ref{pr})   an $h$ is not added to $x$ in any $u^{\ze}(\cd)$ to which $\De_{ x}^{ h}$ is applied. Consequently
  \be
 \De_{ x_{ i}}^{ h}\De_{ y_{ i}}^{ h}u^{\ze,\sharp}( x_{ i}-y_{ i})\Big |_{y_{ i}=x_{ i}}=\De^{h}\De^{-h}u^{\ze}(0).\label{7.9}
 \ee
 
  Continuing,
 we emphasize that the property that  $x_{ i}$ is  a  bound variable depends only on $\si$. The property that  $x_{ i}$ is  a a singular variable depends on the pair $\si,a$. Let
 \begin{equation}
 S(\si,a)=\{i\,|\, \mbox{$x_{ i}$ is a  singular variable for $\si,a$}\}.\label{nm.2}
 \end{equation}

 Consider a  term in  (\ref{nm.1}), with $S(\si,a)=J\subseteq A$. Then for each $i\in J$ we have a unique $k_{i}\in [1,m+|A|]$ such that  $\{\si(k_{i}-1),\si(k_{i})\}=\{x_{i},y_{i}$\}. Let $K=\{k_{i},\,i\in J\}$.      Using (\ref{7.9}),   we see that the contribution of $\si,a$ in  the   second line in (\ref{nm.1})  is:   
\begin{eqnarray} 
 V(\si,a)&:=& \prod_{ j=1}^{ m+|A|}\(\De^{ h}_{ \si(j)}\)^{a_{ 1}(j)}
\(\De^{ h}_{ \si(j-1)}\)^{a_{ 2}(j)}\,u^{\ze,\sharp}(  \si(j)- \si(j-1))\,\,\Big |_{y_{ i}=x_{ i},\,\forall i  }
\nn\\
& =&\(\De^{h}\De^{-h}u^{\ze}(0)\)^{ | J |} \prod_{ j=1,j\notin K}^{ m+|A|}\(\De^{ h}_{ \si(j)}\)^{a_{ 1}(j)}
\(\De^{ h}_{ \si(j-1)}\)^{a_{ 2}(j)}\nn\\
&&\hspace{1.6in}\times \,u^{\ze,\sharp}(  \si(j)- \si(j-1))\Big |_{y_{ i}=x_{ i},\,\forall i }. \label{pr.1}
\end{eqnarray}
 
Let $ \mathcal{I}(\si)$ denote the set of   all
$\si'\in \mathcal{B}_{A}$ which can be obtained from $\si$ by interchanging $\si(k_{i}-1)$ and $\si(k_{i})$ for some   set of  the elements  $i\in J$. 
  Clearly  $V(\si',a)=V(\si,a)$ for all $\si'\in \mathcal{I}(\si)$. Since $|\mathcal{I}(\si)|=2^{|J|}$ we see that the contribution    in  the   second line in (\ref{nm.1}) obtained by   summing over all $\si'\in \mathcal{I}(\si)$ is:   
\begin{eqnarray} 
&&V(\mathcal{I}(\si),a)
=\(2\De^{h}\De^{-h}u^{\ze}(0)\)^{ | J |} \label{pr.1j} \\
&&\hspace{0.4in}\times  \prod_{ j=1,j\notin K}^{m+|A|}\(\De^{ h}_{ \si(j)}\)^{a_{ 1}(j)}
\(\De^{ h}_{ \si(j-1)}\)^{a_{ 2}(j)}  \,u^{\ze,\sharp}(  \si(j)- \si(j-1))\Big |_{y_{ i}=x_{ i},\,\forall i }\nn
\end{eqnarray} 
 
  In what follows  given  $\si \in \mathcal{B}_{ A  }$, we   write it as a vector $(\si(1),\ldots,\si(m+|A|))\in R^{m+|A|}$.  For any $J\subseteq A$ we define $\si_{A-J} \in \mathcal{B}_{ A-J  }$, by deleting the components  $y_{i}$, $ i\in J$ from $(\si(1),\ldots,\si(m+|A|))$. We only use this latter notation when $J$ is contained in the set of singular variables of some $\si,a$.

  As an example of the relationship between $\si$ and $\II(\si)$ let   $m=3$,  $A=\{1,2,3\}$, $\si=(x_{1},x_{2},y_{2},y_{3},x_{3},y_{1})$ and   $J=\{2,3\}$. Then $\mathcal{I}(\si)$
consists of the four bijections
\bea 
  \si= \si_{1}&=&(x_{1},x_{2},y_{2},y_{3},x_{3},y_{1})\label{7.11}\\
   \si_{2}&=&(x_{1},y_{2},x_{2},y_{3},x_{3},y_{1})\nn\\
\si_{3}&=&(x_{1},x_{2},y_{2},x_{3},y_{3},y_{1})\nn\\
\si_{4}&=&(x_{1},y_{2},x_{2},x_{3},y_{3},y_{1})\nn .
 \eea
 Also, in the notation just defined
  $\si_{\{1,3\}}=(x_{1},x_{2} ,y_{3},x_{3},y_{1})$, $\si_{\{1,2\}}=(x_{1},x_{2},y_{2} ,x_{3},y_{1})$ and $\si_{\{1\}}=(x_{1},x_{2},x_{3},y_{1})$. 
 
In the notation just defined, we write (\ref{pr.1j}) as
 \begin{eqnarray} 
V(\mathcal{I}(\si),a)
& =&\(2\De^{h}\De^{-h}u^{\ze}(0)\)^{ | J |}\prod_{ j=1 }^{m+|A-J|}\(\De^{ h}_{ \si_{A-J}(j)}\)^{a'_{ 1}(j)}
\(\De^{ h}_{ \si_{A-J}(j-1)}\)^{a'_{ 2}(j)} \nn \\
&& \hspace{.5 in} \,u^{\ze,\sharp}(  \si_{A-J}(j)- \si_{A-J}(j-1))\Big |_{y_{ i}=x_{ i},\,\forall i }\label{npr.1j}
\end{eqnarray} 
  where $a'$ is obtained from $a=\{(a_{1}(j),a_{2}(j))\}_{j=1}^{m+|A|}$ by deleting from  $a$ the pairs $(a_{1}(j),a_{2}(j))$ for $j\in K$, and renumbering the remaining terms in increasing order.  

Note   that in applying the  product formula for difference operators  (\ref{pr}) we can choose which function plays the role of $f$, and which the role of $g$. When   $x_{i}$ is a bound variable, that is both $x_{ i}, y_{i}$ appear in the same $u^{\ze}(\cd)$, and we apply    (\ref{pr})  to expand  $\De_{ x_{ i}}^{ h}$,  we take  $g$ to be $u^{\ze}( y_{ i}-x_{ i})$. That is, we take
\bea
&&\De_{ x_{ i}}^{ h}u^{\ze}( x_{ i}-a) u^{\ze}( x_{ i}-{y_{i}}) =\\
&&\qquad \De_{x_{1}}^{h} u^{\ze}(x_{i}-a) u^{\ze}( x_{ i}+h-{y_{i}})+ u^{\ze}( x_{ i}-a) \De_{x_{1}}^{h}u^{\ze}( x_{ i} -{y_{i}})\nn,
\eea
 and
similarly when we apply  (\ref{pr})  to expand  $\De_{ y_{ i}}^{ h}$.   Thus if $x_{i}$ is a singular variable and we apply  $\De_{ x_{ i}}^{ h} \De_{ y_{ i}}^{ h}$ by the above rule, and then set $y_{i}=x_{i}$,   the term that contains  $\De^{h}\De^{-h}u^{\ze}(0)$  is
  \be
u^{\ze}( x_{ i}-a)   \De^{h}\De^{-h}u^{\ze}(0)      u^{\ze}( b-x_{ i}).
 \ee 
Note that  there are no $\pm h$ terms added   to the $y_{i}$ or $x_{i}$. Because of this we  see that 
 \begin{eqnarray} 
 & &
\sum_{\{a\,|\,S(\si,a)=J\subseteq A \}}V(\mathcal{I}(\si),a)
=\(2\De^{h}\De^{-h}u^{\ze}(0)\)^{ | J |}  \label{npr.3j} \\
&&  \hspace{1 in}\lc \(\prod_{ i\in A-J}^{\hspace{.2 in}'}\De_{ x_{ i}}^{ h}\De_{ y_{ i}}^{
h}\) \,u^{\ze}(  \si_{A-J}(j)- \si_{A-J}(j-1))\rc \Big |_{y_{ i}=x_{ i},\,\forall i }\nn
\end{eqnarray}
where the notation $\prod^{'}$ indicates that when we  use (\ref{pr}) to expand the second line of (\ref{npr.3j}) we do not apply  
both $\De_{ x_{ i}}^{ h}\De_{ y_{ i}}^{
h}$ to the same factor $u^{\ze}(\cd)$. This is because all the singular variable have been removed from the  $S(\si,a)$.  The significance of this representation is that it does not contain any ambiguous terms $u^{\ze,\sharp}(\cd)$. 

  For $J\subseteq A$, let $\psi\in \mathcal{B}_{A-J}$. We write $\psi$ as a vector in $R^{m+|A-J|}$ whose components consist of a permutation of the $m+|A-J|$ elements   $\{ x_{ i}, y_{ i},  i\in A-J\}\cup \{ x_{ i}, i\in  (A-J)^{ c} \}$.  Let $\si $ be    obtained from this vector by inserting a    component $y_{i}$, following  $x_{i}$,
for each $ i\in J $.   Considering the way $\si_{A-J}$ was obtained from $\si$, (see the paragraph following the one containing (\ref{7.11})), it clear that for this $\si$ we have $ \si_{A-J}=\psi$. It then follows from this and (\ref{npr.3j}) that we can rewrite (\ref{nm.1}) as 
\begin{eqnarray}
&&\sum_{ J\subseteq A  }\(2\De^{h}\De^{-h}u^{\ze}(0)\)^{ | A^{ c} |}\(2\De^{h}\De^{-h}u^{\ze}(0)\)^{ | J |}
\label{npr.4}\\
&&\hspace{.4 in}\lc \(\prod_{ i\in A-J}^{\hspace{.2 in}'}\De_{ x_{ i}}^{ h}\De_{ y_{ i}}^{
h}\)
\sum_{ \si\in \mathcal{B}_{A-J}}\prod_{ j=1}^{m+|A-J|}u^{\ze}( \si ( j)-\si ( j-1))\rc \Big |_{y_{ i}=x_{ i},\,\forall i } \nonumber\\
&&=\sum_{ J\subseteq A  }\(2\De^{h}\De^{-h}u^{\ze}(0)\)^{ | (A-J)^{ c} |} 
\nn\\
&&\hspace{.4 in}\lc \(\prod_{ i\in A-J}^{\hspace{.2 in}'}\De_{ x_{ i}}^{ h}\De_{ y_{ i}}^{
h}\)
\sum_{ \si\in \mathcal{B}_{A-J}}\prod_{ j=1}^{m+|A-J|}u^{\ze}( \si ( j)-\si ( j-1))\rc \Big |_{y_{ i}=x_{ i},\,\forall i }. \nonumber
\end{eqnarray}

Hence by (\ref{1.16g})--(\ref{nm.1}), for any
integer
$m$ we have  
\bea 
&& E\(\( \int ( L^{ x+h}_{ \la_{\ze}}- L^{ x}_{ \la_{\ze}})^{ 2}\,dx-2\De^{h}\De^{-h}u^{\ze}(0) \int  L^{ x}_{
\la_{\ze}}\,dx\)^{ m}
\)\label{npr.5}\\  &&\qquad=\sum_{A\subseteq \{ 1,\ldots,m\}} ( -1)^{|A^{c}|}\sum_{ J\subseteq A  }\int \phi(A-J)\,dx, \nonumber
\eea
where the set function $\phi$ is defined by 
\begin{eqnarray}
&&\phi(D):=\(2\De^{h}\De^{-h}u^{\ze}(0)\)^{ | D^{ c} |} 
\label{npr.6}\\
&&  \hspace{.4 in}\lc \(\prod_{ i\in D}^{\hspace{.2 in}'}\De_{ x_{ i}}^{ h}\De_{ y_{ i}}^{
h}\)
\sum_{ \si\in \mathcal{B}_{D}}\prod_{ j=1}^{m+|D|}u^{\ze}( \si ( j)-\si ( j-1))\rc \Big |_{y_{ i}=x_{ i},\,\forall i }. \nonumber
\end{eqnarray}
It follows from Principle of Inclusion--Exclusion, \cite[p. 66, (8)]{comb}, that
\be
\sum_{A\subseteq \{ 1,\ldots,m\}} ( -1)^{|A^{c}|}\sum_{ J\subseteq A  }\phi(A-J)=\phi(\{ 1,\ldots,m\}).\label{7.21}
\ee

  \medskip	  Referring to   (\ref{npr.5})--(\ref{7.21}) we see that  to estimate   (\ref{1.16g}) we need only consider $A=\{ 1,\ldots,m\}$ and 
those cases  in which each of the $2m$ difference operators $\De^{ h}$  are assigned either to a unique factor  
$u^{\ze}(\cd)$, or if two difference operators are assigned to the same factor $u^{\ze}(\cd)$, it is not of the form   $u^{\ze}(
0)$.  Therefore, we see that 
\bea && E \(\( \int ( L^{ x+h}_{ \la_{\ze}}- L^{ x}_{ \la_{\ze}})^{ 2}\,dx-2\De^{h}\De^{-h}u^{\ze}(0)\int  L^{ x}_{
\la_{\ze}}\,dx\)^{ m}
\)\label{1.20g}\\ &&\qquad=2^{ m}\sum_{ \pi\in \mathcal{D},a}\int \mathcal{T}^{\sharp}_{h}( x;\,\pi,a)\,dx\nn,
\eea where 
\begin{equation}
\mathcal{T}^{\sharp}_{h}( x;\,\pi,a) =\prod_{ j=1}^{ 2m}\(\De^{ h}_{ x_{ \pi( j)}}\)^{a_{ 1}(j)}
\(\De^{ h}_{ x_{ \pi( j-1)}}\)^{a_{ 2}(j)}\,u^{\ze, \sharp}( x_{\pi(j)}- x_{\pi(j-1)})\label{1.21g}
\end{equation}
  and the sum runs over   $\mathcal{D}$, the set of all  maps $\pi\,:\,[1,\ldots, 2m]\mapsto
[1,\ldots, m]$ with $|\pi^{ -1}(i )|=2$ for each $i$, and all  
$a=(a_{ 1},a_{ 2})\,:\,[1,\ldots, 2m]\mapsto \{ 0,1\}\times \{ 0,1\}$ with the
property that for each $i$ there   are  exactly two factors of the form $\De^{
h}_{ x_{i}}$ in (\ref{1.21g}), and if $a( j)=( 1,1)$ for any $j$, then 
$ x_{\pi(j)}\neq x_{\pi(j-1)}$. The factor $2^{ m}$ in (\ref{1.20g})  comes from 
  the fact that $|\pi^{ -1}(i )|=2$ for
each $i$.  

  \medskip

 	It follows from (\ref{7.1}),   (\ref{7.2}), and (\ref{1.20g}) that to obtain (\ref{7.54a}) it suffices to show that
\begin{equation}
   \lim_{h\to 0}h^{-3/2}2^{ m}\sum_{ \pi\in \mathcal{D},a}\int \mathcal{T}^{\sharp}_{h}( x;\,\pi,a)\,dx\label{7.26}
   \end{equation}
   is equal to the right-hand side of (\ref{7.54a}).
 To simplify the proof we first show this with $\mathcal{T}^{\sharp}_{h}( x;\,\pi,a)$ replaced by
\begin{equation}
\mathcal{T}_{h}( x;\,\pi,a) =\prod_{ j=1}^{ 2m}\(\De^{ h}_{ x_{ \pi( j)}}\)^{a_{ 1}(j)}
\(\De^{ h}_{ x_{ \pi( j-1)}}\)^{a_{ 2}(j)}\,u^{\ze}( x_{\pi(j)}- x_{\pi(j-1)}).\label{1.21gbn}
\end{equation}
At the conclusion of this proof we explain why we have the same limits when $\TT_{h}(\cd)$ is replaced by $\TT_{h}^{\sharp}(\cd)$.

\medskip	From this point on the proof is very similar to the   proof  of Lemma \ref{lem-expind}. Let  $m=2n$.
  Consider the multigraph $G_{\pi }$ whose vertices consist of   
$\{1,\ldots, 2n\}$ and  we assign an edge 
between the vertices $\pi (2j-1) $ and $ \pi (2j)$ for each $j=1,\ldots,2n$. Each vertex is connected to two edges, and it is possible to have two edges between any two vertices $i,j$.   Note that the connected components $C_{j}$, $j=1,\ldots, k$ of $G_{\pi}$ consist of cycles.

 \subsection{$a =e$ and all cycles are of order two}\label{ss-3.1t}
 
When $a =e$,  (defined just before Subsection \ref{ss3.1}),  we have
 \begin{equation}\qquad
\mathcal{T}_{h}( x;\,\pi,e) =\prod_{ j=1}^{ 2n}u^{\ze}( x_{\pi(2j-1)}- x_{\pi(2j-2)})\,\De^{ h}\De^{- h}\,u^{\ze}( x_{\pi(2j)}- x_{\pi(2j-1)}).\label{91.1}
\end{equation}
Assume now that, in addition, all cycles are of order two.

 Let $\mathcal{P}=\{(l_{2i-1},l_{2i})\,,\,1\leq i\leq n\}$ be a pairing of the integers $[1,2n]$. Let $\pi\in \mathcal{D}$,   (defined just after (\ref{1.21g})),   be     such that for each $1\leq j\leq 2n$, 
$\{\pi(2j-1), \pi(2j)\}=\{l_{2i-1},l_{2i}\}$ for some, necessarily unique, $ 1\leq i\leq n$. In this case we say that $\pi$ is compatible with the pairing $\mathcal{P}$ and write this  as $ \pi \sim \mathcal{P}$. (Note that  when we write $\{\pi(2j-1), \pi(2j)\}=\{l_{2i-1},l_{2i}\}$ we mean as two  sets, so, according to what $\pi$ is, we may have  $\pi(2j-1)=l_{2i-1}$, $\pi(2j )=l_{2i }$ or $\pi(2j-1)=l_{2i}$, $\pi(2j )=l_{2i-1 }$.)  Whenever $\pi\in \mathcal{D}$   is   such that $G_{\pi}$ consists  only of cycles of order two,     $ \pi \sim \mathcal{P}$, for some pairing $\mathcal{P}$ of the integers $[1,2n]$.
In this case we have 
 \begin{equation} \qquad
\mathcal{T}_{h}( x;\,\pi,e) =\prod_{ i=1}^{ n}\(\De^{ h}\De^{- h}\,u^{\ze}( x_{l_{2i}}- x_{l_{2i-1}})\)^{2}\prod_{ j=1}^{ 2n}u^{\ze}( x_{\pi(2j-1)}- x_{\pi(2j-2)})\,.\label{91.2}
\end{equation}

 Following the proof  of Lemma \ref{lem-expind} we first show that 
\begin{equation}
\int \mathcal{T}_{h}( x;\,\pi,e)\prod_{j=1}^{2n}\,dx_{j} =\int \mathcal{T}_{1,h}( x;\,\pi,a)\prod_{j=1}^{2n}\,dx_{j}+O(h^{3n+1})\label{91.3}
\end{equation}
where
\begin{eqnarray}\qquad
\mathcal{T}_{1,h}( x;\,\pi,e)&=&\prod_{ i=1}^{ n}\(1_{\{|x_{l_{2i}}-x_{l_{2i-1}}|\leq h\}}\)\(\De^{ h}\De^{- h}\,u^{\ze}( x_{l_{2i}}- x_{l_{2i-1}})\)^{2}
\nn\\
&&\qquad   \times \prod_{ j=1}^{ 2n}u^{\ze}( x_{\pi(2j-1)}- x_{\pi(2j-2)}).\label{91.4}
\end{eqnarray}

To prove (\ref{91.3}) we proceed as in (\ref{f9.31})--(\ref{f9.32a}), and   see that it suffices to show that   for $A\subseteq [1,\ldots,n]$ and $|A^{c}|\geq 1$,
\begin{eqnarray} 
&& \int \prod_{i\in A} 1_{\{|x_{l_{2i}}-x_{l_{2i-1}}|\leq h\}}\prod_{i\in A^{c}} 1_{\{|x_{l_{2i}}-x_{l_{2i-1}}|\geq h\}} \(\De^{ h}\De^{- h}\,u^{\ze}( x_{l_{2i}}- x_{l_{2i-1}})\)^{2}
\nn\\
&& \qquad  \times \prod_{ j=1}^{ 2n}u^{\ze}( x_{\pi(2j-1)}- x_{\pi(2j-2)})\prod_{j=1}^{2n}\,dx_{j}=O(h^{3n+1}).\label{7.18}
\end{eqnarray}
 To    show this we first  choose $j_{k}$, $\,k=1,\ldots,n$, so that 
\be \{x_{\pi(2j_{k}-1)}- x_{\pi(2j_{k}-2)},\,k=1,\ldots,n \}\cup \{x_{l_{2i}}- x_{l_{2i-1}},\,i=1,\ldots,n \}\label{7.19}\ee
spans $R^{2n}$.   Let $y_{i}$, $ i=1,\ldots,2n$,  denote the $2n$  variables   in (\ref{7.19}). We make the change of variables in   (\ref{7.18}) to  $\{y_{1},\ldots, y_{2n}\}$.  We then bound those  terms in $u^{\ze}( x_{\pi(2j-1)}- x_{\pi(2j-2)})$, $j=1,\ldots,2n$, that do not map into $u^{\ze}( y_{i} )$, for some $i=1,\ldots,2n$; (see (\ref{pot.1w}).) 
We are then left with an easy integral and using (\ref{1.30g}), and (\ref{1.30gb}) and the fact that $u^{\ze}(\cd)$ is integrable we get (\ref{7.18}).

\medskip	
  Analogous to (\ref{f9.36}) and (\ref{f9.37}) we now study  
\begin{equation} \hspace{.4 in}
\int \mathcal{T}_{1,h}( x;\,\pi,e)\prod_{j=1}^{2n}\,dx_{j}. \label{f91.36}
\end{equation}
 Recall that for  each $1\le j\le 2n$, 
$\{\pi(2j-1), \pi(2j)\}=\{l_{2i-1},l_{2i}\} $, for some $1\le i\le n$.  We identify these relationships by setting $i=\si (j) $ when  $\{\pi(2j-1), \pi(2j)\}=\{l_{2i-1},l_{2i}\} $.  In the present situation,  in which all cycles are of order two,  we have $\si:\,[1,2n]\mapsto [1,n]$,  with $|\si^{-1}(i)|=2$,  for each $1\leq i\leq n$. We write
\begin{eqnarray}
&&\prod_{ j=1}^{2 n}\, u^{\ze}(x_{\pi(2j-1)}-x_{\pi(2j-2)})
 \label{f91.37}\\
&&\qquad=\prod_{ j=1}^{ 2n}\,\( u^{\ze}(x_{l_{2\si (j)-1}}-x_{l_{2\si (j-1)-1}})+\De^{h_{j}}u^{\ze}(x_{l_{2\si (j)-1}}-x_{l_{2\si (j-1)-1}})\) ,  \nn
\end{eqnarray}
where $h_{j}=(x_{\pi(2j-1)}-x_{l_{2\si (j)-1}})+(x_{l_{2\si (j-1)-1}}-x_{\pi(2j-2)})$. 
Note that because of the presence of  the term $\prod_{i=1}^{n}\(1_{\{|x_{l_{2i}}-x_{l_{2i-1}}|\leq h\}}\)$ in the integral in (\ref{f91.36})  we need only be concerned with values of $|h_{j}|\leq 2h$,  $1\le j\le 2n$.

  Following    (\ref{4.28})--(\ref{f9.41z})    we   see that 
\begin{eqnarray}
&&\int \mathcal{T}_{1,h}( x;\,\pi,e)\prod_{j=1}^{2n}\,dx_{j}
\label{91.6}\\
&&=\int   \prod_{ i=1}^{ n}\(\De^{ h}\De^{- h}\,u^{\ze}( x_{l_{2i}}- x_{l_{2i-1}})\)^{2} \prod_{ j=1}^{ 2n}\,  u^{\ze}(x_{l_{2\si (j)-1}}-x_{l_{2\si (j-1)-1}})\prod_{j=1}^{2n}\,dx_{j}\nonumber\\
&&\hspace{3 in}+O(h^{3n+1})\nonumber,
\end{eqnarray}
where $x_{-1}=0$.   

We now estimate the integral in (\ref{91.6}).
 Using translation invariance and then (\ref{1.30g}) we have
\begin{eqnarray}
\lefteqn{\int   \prod_{ i=1}^{ n}\(\De^{ h}\De^{- h}\,u^{\ze}( x_{l_{2i}}- x_{l_{2i-1}})\)^{2} \prod_{ j=1}^{ 2n}\,  u^{\ze}(x_{l_{2\si (j)-1}}-x_{l_{2\si (j-1)-1}})\prod_{j=1}^{2n}\,dx_{j}
\nn}\\
&& =\int   \prod_{ i=1}^{ n}\(\De^{ h}\De^{- h}\,u^{\ze}( x_{l_{2i}})\)^{2} \prod_{ j=1}^{ 2n}\,  u^{\ze}(x_{l_{2\si (j)-1}}-x_{l_{2\si (j-1)-1}})\prod_{k=1}^{2n}\,dx_{l_{k}}  \nonumber\\
&& =( 8/3+O( h))^{n}h^{ 3n}\int     \prod_{ j=1}^{ 2n}\,  u^{\ze}(x_{l_{2\si (j)-1}}-x_{l_{2\si (j-1)-1}})\prod_{k=1}^{n}\,dx_{l_{2k-1}}. \qquad\label{91.7}
\end{eqnarray}

  We set $y_{k}=x_{l_{2k-1}}$ and write the   last line of (\ref{91.7}) as
\be
( 8/3)^{n}h^{ 3n}\int     \prod_{ j=1}^{ 2n}\,  u^{\ze}(y_{\si (j)}-y_{\si (j-1)})\prod_{k=1}^{n}\,dy_{k} +O(h^{3n+1}).\label{7.24}
\ee
  It follows from (\ref{91.3}) and (\ref{91.6})--(\ref{7.24}) that 
\begin{eqnarray}
\lefteqn{\int \mathcal{T}_{h}( x;\,\pi,e)\prod_{j=1}^{2n}\,dx_{j}
\label{91.8}}\\
&& =( 8/3)^{n}h^{ 3n}\int     \prod_{ j=1}^{ 2n}\,  u^{\ze}(y_{\si (j)}-y_{\si (j-1)})\prod_{k=1}^{n}\,dy_{k} +O(h^{3n+1}) \nonumber,
\end{eqnarray}
where $y_{0}=0$.
 
Let $\mathcal{M}$ denote the set of maps $\si$   from $[1,\ldots,2n]$ to $[1,\ldots,n]$ such that $|\si^{ -1}( i)|=2$ for all $i$.
For each pairing $\mathcal{P}$   of $[1,\ldots,2n]$, any 
  $\pi\in \mathcal{D}$ that is compatible with $\PP$, (i.e. 
$\pi\sim \mathcal{P}$),  gives rise    to such a map $\si\in \mathcal{M}$.  Furthermore,   any of the $2^{ 2n}$  maps in $\mathcal{D}$ obtained from $\pi $ by permuting the
$2$ elements in any  of the $2n$ pairs $\{\pi (2j-1), \pi (2j)\}$, give rise to the same map $\si$.   In addition, for any  $\si'\in \mathcal{M}$, we can   reorder  the $2n$ pairs of $\pi$ to obtain a new   $\pi'\sim \mathcal{P}$ which gives rise to $\si'$. Thus we have shown that  
\begin{eqnarray}  \lefteqn{\sum_{ \pi\sim  \mathcal{P}}\int\mathcal{T}_{h}(
x;\,\pi,e)\prod_{ j=1}^{ 2n}\,dx_{j}\label{1.35g}}\\ &&=  \( {32\over 3}h^{ 3}\)^{ n} \sum_{ \si\in  \mathcal{M}}\int     \prod_{ j=1}^{ 2n}\,  u^{\ze}(y_{\si (j)}-y_{\si (j-1)})\prod_{k=1}^{n}\,dy_{k}+O(h^{3n+1})\nn\\ &&  =  \( {16\over 3}h^{ 3}\)^{ n} E\lc\(\int (L^{ x}_{ \la_{\ze}})^{ 2}\,dx\)^{
n}\rc+O(h^{3n+1})\nn
\end{eqnarray}
where the last line follows from Kac's moment formula.   The factor  $2^{-  n}$ that appears in the transition from the   second to the third line in  (\ref{1.35g})  is due to  
  the fact that $|\si^{ -1}(i )|=2$ for
each $i$; (see (\ref{1.20g})).  

  Let $\mathcal{G}_{2}$ denote the set of $\pi\in \mathcal{D}$ such that all cycles of   the graph $G_{\pi}$ have order two. Since every such $\pi$ is compatible with some pairing 
$\mathcal{P}$, and there are $  \displaystyle  {( 2n)! \over
2^{ n}n!}$ such pairings, we see that   \bea&&   \sum_{ \pi \in \mathcal{G}_{2}}\int\mathcal{T}_{h}(
x;\,\pi,e)\prod_{ j=1}^{ 2n}\,dx_{j}\label{1.35gq}\\
&&\qquad = {( 2n)! \over
2^{ n}n!}\( {16\over 3}h^{ 3}\)^{ n} E\lc\(\int (L^{ x}_{ \la_{\ze}})^{ 2}\,dx\)^{
n}\rc+O(h^{3n+1}). \nn
\eea

 \subsection{${\bf a=e }$ and all cycles are not of order two and {${\bf a\neq e }$}}\label{ss-3.2t} 
 
  We follow closely the argument in Subsection \ref{ss-3.2} to show that  
\begin{equation}
\sum_{ \pi \not \in \mathcal{G}_{2}}\Big |\int\mathcal{T}(
x;\,\pi,e)\prod_{ j=1}^{ 2n}\,dx_{j}\Big |=O(h^{3n+1}).\label{91.10}
\end{equation}

   Let the cycles $C_{j}=\{j_{1},\ldots, j_{ l(j)}\}$ of $G_{\pi}$ be written in cyclic order where $l(j)=|C_{j} |$. Note that
$\sum_{j=1}^{k}l(j)=2n$.

Since we    only need an upper bound, we  take absolute values in the integrand   to see that
\begin{eqnarray}
\lefteqn{\bigg|\int \mathcal{T}_{h}( x;\,\pi,e) \,\prod_{ j=1}^{ 2n}\,dx_{j}\bigg|
\nn}\\
&&   \leq   \int    \prod_{j=1}^{k}\( w^{\ze}(x_{j_{2}}-x_{j_{1}})\cdots    w^{\ze}(x_{j_{l(j)}}-x_{j_{l(j)-1}})w^{\ze}(x_{j_{1}}-x_{j_{l(j)}})   \) \nn\\ &&\hspace{
.5in}\times\prod_{ j=1}^{2 n}\, u^{\ze} (x_{\pi(2j-1)}-x_{\pi(2j-2)})   \, \prod_{ j=1}^{ 2n}\,dx_{j} \label{f9.51z},
\end{eqnarray}
where 
  $w^{\ze}(x)$ is defined in (\ref{4.20}).
  Note that we group the functions $w$ according to the cycles.

  We now follow the paragraph containing (\ref{f9.53}) verbatim until the end of Subsection \ref{ss-3.2}, except that we replace $u$ and $w$ by $u^{\ze}$ and $w^{\ze}$, to get (\ref{91.10}).

\medskip	  
  When $a\ne e$  
\begin{equation}
\sum_{\pi }\sum_{  a \neq e}\Big |\int \mathcal{T}_{h}( x;\,\pi ,a ) \,\prod_{ j=1}^{ 2n}\,dx_{j}\Big |=O(h^{3n+1}). \label{91.20}
\end{equation} 
  This follows easily by obvious modifications of the proof in   Subsection \ref{ss3.4}, similar to the modifications of the proof in 
  Subsection \ref{ss-3.2} that gives   (\ref{91.10}).

 \medskip	  We now note that it follows from the arguments in the final three paragraphs of the proof of  Lemma \ref{lem-expind}, on page \pageref{page24}, that for $m$ even we obtain the same 
  asymptotic behavior  when we replace   $\mathcal{T}_{h}( x;\,\pi,a)$ by $\mathcal{T}^{\sharp}_{h}( x;\,\pi,a)$, and also, that we get the right-hand side of (\ref{7.54a}) for odd moments.
  
\medskip	   Summing up, we have shown that the only non-zero limits in (\ref{7.26}) come from (\ref{1.35gq}) when $m$ is even. Using this in (\ref{7.26}), in which we multiply by $2^{2n}$,  we see that (\ref{7.26}) is equal to the right-hand side of (\ref{7.54a}). \qed

\section{Expectation}\label{sec-8}

\begin{lemma}\label{lem-vare} For $h\geq 0$  
\begin{equation}
E\(\int ( L^{ x+h}_{1}- L^{ x}_{ 1})^{ 2}\,dx\)=4h+O(h^{2}),\label{expe}
\end{equation}
as $h\to 0$.
Equivalently,  
\begin{equation}
E\(\int ( L^{ x+1}_{t}- L^{ x}_{ t})^{ 2}\,dx\)=4t+O(t^{1/2}),\label{expet}
\end{equation}
as $t\to\ff.$
\end{lemma} 

 \Proof   By the     Kac  moment formula 
\begin{eqnarray}
&&E\(\int ( L^{ x+h}_{1}- L^{ x}_{ 1})^{ 2}\,dx\)
\label{kacv.3}\\
&&\qquad =2 \int \int_{\{\sum_{i=1}^{2}r_{i}\leq 1\}}\De^hp_{r_{1}}(x)\De^h p_{r_{2}}( 0)\,dr_{1}\,dr_{2}\,dx \nonumber\\
&&\quad \qquad+2 \int \int_{\{\sum_{i=1}^{2}r_{i}\leq 1\}}p_{r_{1}}(x)\De^h \De^{-h} p_{r_{2}}( 0)\,dr_{1}\,dr_{2}\,dx. \nonumber
\end{eqnarray}
When we  integrate with respect to  $x$ we get zero in the first integral and one in the second. Consequently
 \begin{eqnarray}
E\(\int ( L^{ x+h}_{1}- L^{ x}_{ 1})^{ 2}\,dx\) &=& 2   \int_{\{\sum_{i=1}^{2}r_{i}\leq 1\}} \De^h \De^{-h} p_{r_{2}}( 0)\,dr_{1}\,dr_{2}\label{kacv.4}\\
 &=& 4  \int_{0}^{1} (1-r)\(p_{r }( 0)-p_{r }( h)\) \,dr \nn.
\end{eqnarray}

Since 
\begin{equation}
\int_{0}^{1} r {1-e^{-h^{2}/2r} \over \sqrt{r}} \,dr\leq \int_{0}^{1} r {  h^{2}/2r\over \sqrt{r}} \,dr=O(h^{2}).\label{kacv.6}
\end{equation}
and
\begin{equation}
\int_{1}^{\ff}  {1-e^{-h^{2}/2r} \over \sqrt{r}} \,dr\leq \int_{1}^{\ff}  {  h^{2}/2r\over \sqrt{r}} \,dr=O(h^{2})\label{kacv.7aa}
\end{equation}
we see that to prove (\ref{expe}) it suffices to show that 
\begin{equation}
    \int_{0}^{\ff} \(p_{r }( 0)-p_{r }( h)\) \,dr=h+O(h^{2}).
   \end{equation}
   This follows from (\ref{pot.1w}) since 
\begin{equation}
    \int_{0}^{\ff} \(p_{r }( 0)-p_{r }( h)\) \,dr=\lim_{\al\to 0}  \int_{0}^{\ff}e^{-\al r}\(p_{r }( 0)-p_{r }( h)\) \,dr.   \end{equation}
Thus we get (\ref{expe});
 (\ref{expet})   follows from the scaling property,   (\ref{scale}).\qed

\def\noopsort#1{} \def\printfirst#1#2{#1}
\def\singleletter#1{#1}
      \def\switchargs#1#2{#2#1}
\def\bibsameauth{\leavevmode\vrule height .1ex
      depth 0pt width 2.3em\relax\,}
\makeatletter
\renewcommand{\@biblabel}[1]{\hfill#1.}\makeatother
\newcommand{\bysame}{\leavevmode\hbox to3em{\hrulefill}\,}

\end{document}